\documentclass[10pt]{article}
\usepackage{latexsym,amsmath,amssymb,graphics,amscd}
\usepackage{hyperref}
\hypersetup{colorlinks=true, citecolor=blue}
\usepackage{stackengine}
\usepackage{graphicx}
\usepackage{subcaption}

\addtocounter{footnote}{1}
\makeatletter
\@addtoreset{table}{bsection}
\def\thetable{\thesection.\@arabic\c@table}
\def\fps@table{h, t}
\@addtoreset{equation}{section}

\makeatother

\newtheorem{theorem}{Theorem}[section]

\newtheorem{lemma}[theorem]{Lemma}

\newtheorem{proposition}[theorem]{Proposition}

\newtheorem{corollary}[theorem]{Corollary}

\newcommand{\bfi}{\bfseries\itshape}

\newsavebox{\savepar}

\makeatletter

\begin{document}

\title{\textbf{Learning strange attractors with reservoir systems}}
\author{Lyudmila Grigoryeva$^{1}$, Allen Hart$^{2}$, and Juan-Pablo Ortega$^{3}$}
\date{}

\maketitle

\makeatletter
\addtocounter{footnote}{1} \footnotetext{%
Department of Statistics, University of Warwick, Coventry CV4 7AL , UK. {\texttt{} }}
\makeatother
\makeatletter
\addtocounter{footnote}{1} \footnotetext{%
Department of Mathematical Sciences, University of Bath, Bath BA2 7AY, UK.
{\texttt{A.Hart@bath.ac.uk} }}
\makeatother
\makeatletter
\addtocounter{footnote}{1} \footnotetext{%
Division of Mathematical Sciences, 
Nanyang Technological University,
21 Nanyang Link,
Singapore 637371.
{\texttt{Juan-Pablo.Ortega@ntu.edu.sg}}}
\makeatother

\begin{abstract}
This paper shows that the celebrated Embedding Theorem of Takens is a particular case of a much more general statement according to which, randomly generated linear state-space representations of generic observations of an invertible dynamical system carry in their wake an embedding of the phase space dynamics into the chosen Euclidean state space. This embedding coincides with a natural generalized synchronization that arises in this setup and that yields a topological conjugacy between the state-space dynamics driven by the generic observations of the dynamical system and the dynamical system itself. This result provides additional tools for the representation, learning, and analysis of chaotic attractors and sheds additional light on the reservoir computing phenomenon that appears in the context of recurrent neural networks.
\end{abstract}


\section{Introduction}
\label{Introduction}

Takens' Theorem \cite{takensembedding} and the associated method of delays have been used and studied for decades as they are powerful tools in the reconstruction of qualitative features of a dynamical system out of time series of low dimensional observations. This result is also at the origin of the development of powerful forecasting tools \cite{sauer1991embedology, kantz:Sreiber}. 

In order to put these results in context and to better motivate the contributions in this paper, we start by recalling Huke's formulation \cite{huke:2006} of Takens' Theorem.

\begin{theorem}[Takens]
Let $M$ be a compact manifold of dimension $q \in \mathbb{N}$ and let $\phi \in {\rm Diff} ^2(M) $ be a twice-differentiable diffeomorphism that satisfies the following two properties:
\begin{description}
\item [(i)] $\phi $ has only finitely many periodic points with periods less than or equal to $2q$.
\item [(ii)] If $m \in M $  is any periodic point of $\phi $ with period $k<2q $, then the eigenvalues of the linear map $T _m \phi ^k: T _m M \longrightarrow T _mM  $ are distinct.
\end{description}
Then for any generic scalar observation function $\omega \in C ^2(M , \mathbb{R} )$, the $(2q+1) $-delay map $\Phi_{(\phi, \omega)}: M \longrightarrow \mathbb{R}^{2q+1}$ defined by
\begin{equation}
\label{delay embedding expression}
\Phi_{(\phi, \omega)}(m):= \left(
\omega(m), \omega \circ \phi (m), \omega \circ \phi^2 (m), \ldots, \omega \circ \phi^{2q} (m)
\right)
\end{equation}
is an embedding in $C ^1(M,\mathbb{R}^{2q+1}) $.
\end{theorem}

The first  consequence of this result is that, since the map $\Phi_{(\phi, \omega)} $ is an embedding, then it is necessarily injective and hence it can be used to represent in $\Phi_{(\phi, \omega)} (M) \subset\mathbb{R}^{2q+ 1 }$ the dynamics induced by $\phi$ on $M$ via the differentiable map $\varphi_{(\phi, \omega)}:= \Phi_{(\phi, \omega)}\circ \phi \circ \Phi_{(\phi, \omega)} ^{-1}:\Phi_{(\phi, \omega)} (M)\subset \mathbb{R}^{2q+ 1 } \longrightarrow \Phi_{(\phi, \omega)} (M)\subset\mathbb{R}^{2q+ 1 } $ (we recall that the inverse function theorem guarantees that the map $ \Phi_{(\phi, \omega)} ^{-1}:\Phi_{(\phi, \omega)} (M) \longrightarrow M$ is differentiable). In view of the expression \eqref{delay embedding expression}, this map takes necessarily the form $\varphi_{(\phi, \omega)}(z _1, \ldots, z_{2q+1}) :=(z _2, z _3, \ldots , h(z _1, \ldots, z_{2q+1}))$, for some differentiable map $h:\Phi_{(\phi, \omega)} (M)\subset\mathbb{R}^{2q+ 1 } \longrightarrow \mathbb{R} $. In this situation, we say that the dynamical systems $(M,\phi) $ and $(\Phi_{(\phi, \omega)} (M),\varphi_{(\phi, \omega)} ) $ are topologically conjugate by the map $\Phi_{(\phi, \omega)}  $. The importance of this representation is that the two systems $(M,\phi) $ and $(\Phi_{(\phi, \omega)} (M),\varphi_{(\phi, \omega)} ) $ have the same $C ^1 $ invariants like Lyapunov exponents, eigenvalues of linearizations, or dimensions of attractors and their computation may be more efficiently carried out in  $\Phi_{(\phi, \omega)} (M)\subset\mathbb{R}^{2q+ 1 } $.

More recently, the remarkable success of recurrent neural networks and reservoir computing \cite{lukosevicius, tanaka:review} in the learning, forecasting \cite{Jaeger04, pathak:chaos, Pathak:PRL, Ott2018}, and classification \cite{carroll2018using} of chaotic attractors of complex nonlinear high-dimensional dynamical systems strongly suggests that these machine learning paradigms have Takens embedding-type properties. This fact has been rigorously established in \cite{hart:ESNs, allen:tikhonov} where the so called {\it Echo State Networks} (ESNs) \cite{Matthews:thesis, Matthews1993, Jaeger04, RC7, RC8, RC20} driven by one-dimensional observations of a given dynamical system on a compact manifold have been shown, under certain hypotheses, to produce dynamics that are topologically conjugate to that of the original system. 

A concept that unifies the recurrent networks and the Takens approaches to the representation of dynamical systems is that of {\bfi  generalized sychronization} (GS), as introduced in \cite{rulkov1995generalized} (see \cite{pecora:synch, ott2002chaos, boccaletti:reports:2002, eroglu2017synchronisation} for self-contained presentations and many references). Generalized synchronizations represent dynamical systems in the space of states $\mathbb{R}^N $ of a state-space map $F: \mathbb{R} ^N\times \mathbb{R} ^d \longrightarrow   \mathbb{R} ^N $, $N, d \in \mathbb{N}$. More specifically, let $(M, \phi)$ be the same dynamical systems as above, with  $M$ compact and $\phi\in {\rm Diff}^1(M) $. Let $\omega \in C ^1(M, \mathbb{R} ^d) $, $d \in \mathbb{N} $,  be a map that encodes $d$-dimensional observations of the dynamical system and define the $(\phi, \omega)$-{\it delay map} $S _{(\phi, \omega)}:M \longrightarrow \ell^{\infty}(\mathbb{R}^d) $ as $S _{(\phi, \omega)}(m):=\left\{\omega(\phi ^t (m))\right\}_{t \in \Bbb Z} $. Consider now the  drive-response system associated to the $\omega $-observations of $\phi$ and determined by the recursions:
\begin{equation}
\label{drive-response system}
\mathbf{x} _t=F\left(\mathbf{x} _{t-1}, S _{(\phi, \omega)}(m)_t\right), \quad \mbox{$t \in \Bbb Z,\, m \in M.$}
\end{equation}
We say that a generalized synchronization occurs in this configuration when there exists a map $f_{(\phi, \omega, F)}:M \longrightarrow \mathbb{R}^N $ such that for any $\mathbf{x} _t  $, $t \in \Bbb Z $, as in \eqref{drive-response system}, it holds that
\begin{equation}
\label{generalized synchronization condition}
 \mathbf{x} _t = f_{(\phi, \omega, F)} (\phi ^t(m)),
\end{equation}
that is, the time evolution of the dynamical system in phase space (not just its observations) drives the response in \eqref{drive-response system}. We emphasize that the definition \eqref{generalized synchronization condition} presupposes that the recursions \eqref{drive-response system} have a (unique) solution, that is, that there exists a sequence $\mathbf{x} \in \ell^{\infty}(\mathbb{R}^N) $ such that \eqref{drive-response system} holds true. When that existence property holds and, additionally, the solution sequence $\mathbf{x}  $ is unique, we say that $F$ has the $(\phi, \omega) $-{\it Echo State Property} (ESP) (see \cite{jaeger2001, Manjunath:Jaeger, manjunath:prsl} for in-depth descriptions of this property). Moreover, in the presence of the $(\phi, \omega) $-ESP, the state map $F$ determines a unique causal and time-invariant filter $U ^F: S _{(\phi, \omega)}(M) \longrightarrow (\mathbb{R} ^N)^{\Bbb Z} $ that associates to each orbit $S _{(\phi, \omega)}(m) $ the unique solution sequence $\mathbf{x} \in (\mathbb{R} ^N)^{\Bbb Z} $ of \eqref{drive-response system}. The existence, continuity, and differentiability of GSs has been established in \cite{RC18} for a rich class of systems that exhibit the so-called fading memory property and that are generated by locally state-contracting maps $F$.

The relevance of these concepts in relation to the embedding of dynamical systems lays in the fact that Takens' Theorem can be easily reformulated in the language of generalized synchronizations. Indeed, we first note that the map 
$\Phi_{(\phi, \omega)}$ introduced in \eqref{delay embedding expression} is the GS  corresponding to the linear state map $F(\mathbf{x}, z):=A \mathbf{x}+ \mathbf{C} z $, with $A$ the lower shift matrix in dimension $2q+1  $ and $\mathbf{C}= (1,0, \ldots,0) ^{\top} \in \mathbb{R}^{2q+1} $. Takens' Theorem can now be stated by saying the the GS $\Phi_{(\phi, \omega)}$ is an embedding for any generic scalar observation function $\omega \in C ^2(M , \mathbb{R} )$. 

{\it The main result in this paper shows that Takens' Theorem  is a particular case of a more general statement that ensures that the GSs associated to {\it generic randomly generated} linear state-space systems of the type $F(\mathbf{x}, z):=A \mathbf{x}+ \mathbf{C} z $, with $A \in \mathbb{M}_{N,N} $, $\mathbf{C} \in \mathbb{R}^N$, and $N \geq 2q+1 $, and driven by generic  observations $\omega \in C ^2(M , \mathbb{R} )$   are embeddings.}

The term {\it generic} is used in the previous statement with two different meanings. First, when we talk about {\it generic} randomly generated linear state-space systems, we mean that the embedding condition holds almost surely when $A \in \mathbb{M}_{N,N} $ and $\mathbf{C} \in \mathbb{R}^N$ are randomly drawn with respect to some probability distribution in a subset of the spaces where those elements are defined. Second, when we write  {\it generic}  observations $\omega \in C ^2(M , \mathbb{R} )$, we mean that they belong to an open and dense subset of $C ^2(M , \mathbb{R} )$ with respect to a Banach topology in that space that we define later on in the paper.

An important consequence of this result is that it sheds light on the good performance of {\bfi  reservoir computing} (RC)  \cite{Jaeger04, lukosevicius, tanaka:review} in the forecasting of dynamical systems. We recall that RC (also found in the literature under other denominations like {\it   Liquid State Machines}~\cite{Maass2000, maass1, Natschlager:117806, corticalMaass, MaassUniversality}) capitalizes on the  idea that there are randomly generated systems that attain universal approximation properties without the need to estimate all their parameters. 
RC has shown unprecedented abilities in the learning of the attractors of complex nonlinear infinite dimensional dynamical systems \cite{Jaeger04, pathak:chaos, Pathak:PRL, Ott2018} and has given rise to forecasting techniques that outperform standard Takens-based strategies. 

Our results explicitly contribute in relation to the  RC phenomenon by showing that the dynamics of generic observations of invertible dynamical systems is almost surely learnable using randomly generated linear reservoir systems with nonlinear readouts (unlike what is common practice in the RC literature, where readouts are linear). Indeed, let $f_{(\phi, \omega, F)}:M \longrightarrow \mathbb{R}^N$ be a GS associated to a randomly generated linear state-space system that, as above, is driven by  generic scalar observations  $\omega \in C ^2(M , \mathbb{R} )$ of $ \phi \in  {\rm Diff} ^2(M)$. Since our results show that $f_{(\phi, \omega, F)} $ is an embedding, it then has an inverse and we can hence construct the readout $h:= \omega \circ \phi \circ f _{(\phi, \omega, F)}^{-1}: f(M) \subset  \mathbb{R}^N \longrightarrow \mathbb{R} $ that, applied to the states $\mathbf{x}_t $ determined by \eqref{generalized synchronization condition} fully characterize the dynamics of the $\omega $-observations $\left\{\omega(\phi ^t (m))\right\}_{t \in \Bbb Z} $ of $\phi$ because $ h(\mathbf{x} _t)= \omega\circ  \phi(f_{(\phi, \omega, F)} ^{-1}(\mathbf{x} _t))= \omega(\phi ^{t+1} (m))$. This observation implies that this dynamics can captured via the  learning of the function $h:= \omega \circ \phi \circ f_{(\phi, \omega, F)} ^{-1} $. This is what we call {\it learnability}  (see, for instance, \cite{lu:bassett:2020, Verzelli2020b, gauthier2021next}).  We emphasize that the regularity properties of the map $f_{(\phi, \omega, F)} $ that we establish later on in the paper guarantee that the readout $h:= \omega \circ \phi \circ f _{(\phi, \omega, F)}^{-1}$ can be efficiently approximated by a universal family (for instance neural networks or polynomials) and explains the good performance of this methodology in the applications cited above. 

The paper is organized as follows. Section \ref{Definitions and preliminary discussion} contains a first introduction to the connection between generalized synchronizations and embeddings and provides existence and regularity statements in the linear case (mainly Proposition \ref{generalized synch with spectral radius}) that are used later on in the paper. Section \ref{Immersive generalized synchronizations} introduces and proves Theorem \ref{Theorem immersion}, which establishes sufficient conditions for a linear system to yield immersive generalized synchronizations for generic observation maps. In Section \ref{Linear reservoir embeddings} we show first (Theorem \ref{the immersions are embeddings}) that basically without additional hypotheses, the globally immersive generalized synchronizations whose existence was proved in Theorem \ref{Theorem immersion} are injective and hence are necessarily embeddings due to the compactness of $M$. Finally, it is also shown (Theorem \ref{Theorem Linear reservoir embeddings}) in this section that randomly generated linear systems (linear reservoirs) yield  synchronization maps $f_{(\phi, \omega,F)} \in C ^2(M, \mathbb{R} ^N ) $ that are almost surely embeddings and are hence amenable to learnability from data. Section \ref{Numerical illustrations} contains a series of numerical illustrations that show the pertinence of the proposed results for attractor reconstruction, filtering in the presence of noise, and forecasting.

\section{Definitions and preliminary discussion}
\label{Definitions and preliminary discussion}

All along this paper we consider an invertible and discrete-time dynamical system determined by a map  $\phi$ that belongs to the set of diffeomorphisms  ${\rm Diff} ^1(M) $ of a finite-dimensional compact manifold $M $. Since later on we need to ensure that $M$ can be endowed with a Riemannian metric $g$, we additionally assume that $M$ is connected, Hausdorff, and second-countable (see \cite[Proposition 2.10]{do:carmo:1993}). The $d$-dimensional observations of the dynamical system are realized by maps $\omega $ that belong, most of the time, to $C ^1(M, \mathbb{R}^d) $. The symbol $TM $ denotes the tangent bundle of  $M$, $T \phi:TM \longrightarrow TM $ the tangent map of $\phi $,  and $D \omega: TM \longrightarrow \mathbb{R}^d   $  the differential of the observation map $\omega $.
Now, for any $f \in C ^1(M, \mathbb{R} ^N)$, define
\begin{equation*}
\left\|Df\right\|_{\infty}=\sup_{m \in M} \left\{ \left\|Df(m)\right\|\right\} \quad \mbox{with} \quad
\left\|Df(m)\right\|=\sup_{\stackanchor{$\scriptstyle \mathbf{v} \in T _mM$}{$\scriptstyle \mathbf{v}\neq {\bf 0}$} } 
\left\{\frac{\left\|Df(m) \cdot \mathbf{v}\right\|}{\left(g(m)(\mathbf{v},\mathbf{v})\right)^{1/2}}\right\}.
\end{equation*}
Analogously, if $\phi:M \rightarrow M $ is a $C ^1 $ map, we can define:
\begin{equation*}
\left\|T \phi\right\|_{\infty}=\sup_{m \in M} \left\{ \left\|T _m \phi\right\|\right\} \quad \mbox{with} \quad
\left\|T _m \phi\right\|=\sup_{\stackanchor{$\scriptstyle \mathbf{v} \in T _mM$}{$\scriptstyle \mathbf{v}\neq {\bf 0}$} } 
\left\{\frac{\left(g(\phi(m))(T _m\phi \cdot \mathbf{v},T _m\phi \cdot \mathbf{v})\right)^{1/2}}{\left(g(m)(\mathbf{v},\mathbf{v})\right)^{1/2}}\right\}.
\end{equation*}

It can be proved by using the results in Chapter 2 of \cite{abraham:robbin} that the norm $\left\|\cdot \right\|_{C ^1}$ defined by
\begin{equation}
\label{definition c1d}
\left\|f \right\|_{C ^1}:= \left\|f\right\|_{\infty}+ \left\|Df\right\|_{\infty}
\end{equation}
endows $C ^1(M, \mathbb{R} ^N)$ with a Banach space structure. Additionally, (see \cite[Theorem 11.2 {\bf (ii)}]{abraham:robbin})  this norm generates a topology in $C ^1(M, \mathbb{R} ^N)$ that is independent of the choice of Riemannian metric $g$ and coincides with the weak and strong topologies introduced in Chapter~2 of \cite{Hirsch:book}.  These notions can be extended to higher order differentiable maps in a straightforward manner.

The embeddings that are at the core of this paper will be constructed using generalized synchronizations associated to linear systems. That is why we start by recalling a result proved in \cite{RC18} in relation with the existence of these objects in a rich variety of situations. The statement requires the following constants defined with respect to a subset $V \subset \mathbb{R}^N $:
\begin{equation}
\label{Ls for F}
\begin{array}{lllrrr}
L_{F _x}&:=&\sup_{(\mathbf{x}, {\bf z})\in V \times \omega(M)} \left\{ \left\|D _x F(\mathbf{x}, {\bf z})\right\|\right\},&L_{F _z}&:=&\sup_{(\mathbf{x}, {\bf z})\in V \times \omega(M)} \left\{ \left\|D _z F(\mathbf{x}, {\bf z})\right\|\right\},\\
L_{F _{xx}}&:=&\sup_{(\mathbf{x}, {\bf z})\in V \times \omega(M)} \left\{ \left\|D _{xx} F(\mathbf{x}, {\bf z})\right\|\right\},&L_{F _{xz}}&:=&\sup_{(\mathbf{x}, {\bf z})\in V \times \omega(M)} \left\{ \left\|D _{xz} F(\mathbf{x}, {\bf z})\right\|\right\},
\end{array}
\end{equation}

\begin{theorem}[Existence and uniqueness of differentiable generalized synchronizations]
\label{differentiable generalized synchronizations theorem}
Let $\phi \in {\rm Diff}^1(M)$ be a dynamical system on the compact manifold $M$ and consider the observation $\omega \in C ^1(M, \mathbb{R}^d) $  and state  $F\in C ^2(D_N \times D_d, D_N) $  maps, with $D _N\subset \mathbb{R} ^N  $ and $D _d \subset {\Bbb R}^d $ open subsets such that $\omega(M) \subset D_d$. Let $V\subset D_N$ be a closed convex subset and suppose that  $F(V \times \omega(M) )\subset V$. Suppose that the bounds for the partial derivatives of $F$  introduced in \eqref{Ls for F} are all finite and that, additionally,
\begin{equation}
\label{condition for lfx in differentiable}
L_{F _x}< \min \left\{1, 1/ \left\|T \phi ^{-1}\right\|_{\infty}\right\}.
\end{equation}
Then there exists a compact and convex subset $W \subset V$ such that  $F(W \times \omega(M))\subset W $ and: 
\begin{description}
\item [(i)] The system determined by $F:W \times \omega(M) \longrightarrow W$ and driven by the $\omega $-observations of $\phi $ has the $(\phi, \omega) $-ESP and  a generalized synchronization $f_{(\phi, \omega,F)}: M \longrightarrow W $ exists and  is well-defined by the relation $U^F(S _{(\phi, \omega)}(m)) _t=f_{(\phi, \omega,F)} \left(\phi ^t (m)\right)$, for any $t \in \Bbb Z,\, m \in M$.
\item [(ii)] The map $f_{(\phi, \omega,F)} $ belongs to $C ^1(M, W) $ and it is the only one that satisfies the identity:
\begin{equation*}
f_{(\phi, \omega,F)} (m)=F \left(f_{(\phi, \omega,F)} (\phi ^{-1}(m)), \omega (m)\right), \quad \mbox{for all} \quad m \in  M.
\end{equation*}
\end{description}
\end{theorem}

If we now consider the linear system 
\begin{equation}
\label{linear system definition}
F(\mathbf{x}, z):=A \mathbf{x}+ \mathbf{C} z, \quad \mbox{with $A \in \mathbb{M}_{N,N} $, $\mathbf{C} \in \mathbb{R}^N$, $N \in \mathbb{N}  $,}
\end{equation}
 in the context of the previous theorem, we obtain the following corollary that is a straightforward consequence of the fact that, in this case, $A= D _xF(\mathbf{x}, z) $, for all $ \mathbf{x} \in \mathbb{R} ^N  $  and $z \in \mathbb{R} $, and hence $L_{F _x}= \left\|A\right\|$. We shall refer to $A$ as the {\it connectivity matrix} and to the vector $\mathbf{C} $ as the {\it input mask}.

\begin{corollary}
\label{GS linear case}
Let $\phi \in {\rm Diff}^1(M)$ be a dynamical system on the compact manifold $M$ and consider the observation map $\omega \in C ^1(M, \mathbb{R}) $. Let $F: \mathbb{R}^N \times \mathbb{R} \longrightarrow \mathbb{R}^N  $ be the linear state map  given by $F(\mathbf{x}, z):=A \mathbf{x}+ \mathbf{C} z $ with  $A \in \mathbb{M}_{N,N} $, $\mathbf{C} \in \mathbb{R}^N$, $N \in \mathbb{N}  $, such that 
\begin{equation}
\label{condition for GS linear}
\left\|A\right\|< \min \left\{1, 1/ \left\|T \phi ^{-1}\right\|_{\infty}\right\}.
\end{equation}
\begin{description}
\item [(i)] The system determined by $F:\mathbb{R}^N \times \omega(M) \longrightarrow \mathbb{R}^N $ and driven by the $\omega $-observations of $\phi $ has the $(\phi, \omega) $-ESP and  a generalized synchronization $f_{(\phi, \omega,F)}: M \longrightarrow \mathbb{R}^N$ given by 
\begin{equation}
\label{functional form GS linear case}
f_{(\phi, \omega,F)}(m)=\sum _{j=0}^{\infty}A ^j \mathbf{C} \omega(\phi^{-j}(m)).
\end{equation}
\item [(ii)] The map $f_{(\phi, \omega,F)} $ belongs to $C ^1(M, \mathbb{R}^N) $ and it is the only one that satisfies the identity:
\begin{equation*}
f_{(\phi, \omega,F)} (m)=A f_{(\phi, \omega,F)} (\phi ^{-1}(m))+ \mathbf{C}\omega (m), \quad \mbox{for all} \quad m \in  M.
\end{equation*}
\end{description}
\end{corollary}

The features of the linear case allow us to prove the existence of generalized synchronizations in situations that go beyond those spelled out in Theorem \ref{differentiable generalized synchronizations theorem} and Corollary \ref{GS linear case}. More specifically, an argument similar to what can be found in Proposition 4.2 in \cite{RC16} allows us to drop the compactness condition on the manifold $M$ and to replace the hypotheses on the design matrix $A$ by more general ones based on its spectral radius $\rho(A)$.

\begin{proposition}
\label{generalized synch with spectral radius}
Let $\phi \in {\rm Diff}^1(M)$ be a dynamical system on the manifold $M$ (not necessarily compact) and consider the observation map $\omega \in C ^1(M, \mathbb{R}) $. Let $F: \mathbb{R}^N \times \mathbb{R} \longrightarrow \mathbb{R}^N  $ be a linear state map  given by $F(\mathbf{x}, z):=A \mathbf{x}+ \mathbf{C} z $ with  $A \in \mathbb{M}_{N,N} $, $\mathbf{C} \in \mathbb{R}^N$, $N \in \mathbb{N}  $.
\begin{description}
\item [(i)] If the spectral radius of $A$ satisfies that $\rho(A)<1$ and $\omega$ maps into a bounded set of $\mathbb{R} $ then the GS $f_{(\phi, \omega,F)}: M \longrightarrow \mathbb{R}^N$ introduced in \eqref{functional form GS linear case} exists and it is a continuous map.
\item [(ii)] Additionally, let $r \in \mathbb{N}$ and suppose that $\phi \in {\rm Diff}^r(M)$ and that there exist constants $k _1, \ldots, k _r \in \mathbb{N}  $ such that $\left\|A^{k _i}\right\| \left\|T ^i \phi ^{-k _i}\right\|_{\infty}<1   $, $\left\|T^i \phi ^{-1}\right\|_{\infty}< \infty$ for all $i \in  \left\{1, \ldots, r\right\}$. Then for any $\omega \in C ^r(M, \mathbb{R}) $ such that  $\left\|D^i\omega\right\|_{\infty}< \infty $, for all $i \in  \left\{1, \ldots, r\right\}$, the map $f_{(\phi, \omega,F)} $ belongs to $ C ^r(M, \mathbb{R}^N) $ and the higher order derivatives are given by:
\begin{equation}
\label{higher order derivatives of GS}
D^if_{(\phi, \omega,F)}(m)=\sum _{j=0}^{\infty}A ^j \mathbf{C} D^i \left(\omega \circ \phi^{-j}\right) (m),  \mbox{ for all $i \in  \left\{1, \ldots, r\right\}$.}
\end{equation}
\item [(iii)] Suppose now that $M$ is compact. In the hypotheses of points {\bf (i)} and {\bf (ii)} above, the map
\begin{equation}
\label{functor gs}
\begin{array}{cccc}
\Theta_{(\phi, F)}: &C ^r(M, \mathbb{R}) & \longrightarrow & C ^r(M, \mathbb{R}^N)\\
	& \omega &\longmapsto & f_{(\phi, \omega,F)}
\end{array}
\end{equation}
is continuous. Moreover, the subsets $\Omega _i$ and $\Omega_e $ of $C ^r(M, \mathbb{R}) $ for which the corresponding GS are immersions and embeddings, respectively, are open.
\end{description}
\end{proposition} 

\noindent\textbf{Proof.\ \ (i)} This statement is obtained out of a combination Weierstrass M-test (see \cite[Theorem 9.6]{Apostol:analysis}) and Gelfand's formula  for the spectral radius  (see \cite{lax:functional:analysis}), that is, 
$\lim\limits_{k \rightarrow \infty}\left\|{A ^k}\right\|^{1/k}=\rho(A)$. Since by hypothesis $\rho(A)<1 $, we can guarantee the existence of a number $k _0 \in \mathbb{N}$ such that $\left\|A ^{k _0}\right\|<1 $, for all $k\geq k _0 $. Consider now the series $\sum _{j=0}^{\infty}A ^j \mathbf{C} \omega(\phi^{-j}(m)) $ in \eqref{functional form GS linear case} that defines $f_{(\phi, \omega,F)}(m)  $. Given that for any $j \in \mathbb{N}  $ there exist $l \in \mathbb{N}  $  and $i \in \left\{0, \ldots, k _0-1\right\} $ such that $A ^j=A^{l k _0+i}$, we then have that,
\begin{equation}
\label{sum for M}
\left\|A ^j \mathbf{C} \omega(\phi^{-j}(m))\right\|\leq \left\|A^{k _0}\right\| ^l \left\|\mathbf{C}\right\|K _A K _\omega,
\end{equation}
with $K _A=\max \left\{1, \left\|A\right\|, \ldots, \left\|A^{k _0 -1}\right\| \right\}$ and $K _\omega \in \mathbb{R}  $ a constant that satisfies that $| \omega (m)|\leq K _\omega  $ for any $m \in M $ and that is available by the boundedness hypothesis on $\omega (M) $.

The inequality \eqref{sum for M} and the Weierstrass M-test guarantee that the series $\sum _{j=0}^{\infty}A ^j \mathbf{C} \omega(\phi^{-j}(m)) $ converges absolutely and uniformly on $M$ and that
\begin{equation*}
\left\|f_{(\phi, \omega,F)} (m)\right\|= \left\|\sum _{j=0}^{\infty}A ^j \mathbf{C} \omega(\phi^{-j}(m))\right\|\leq 
\sum _{l=0}^{\infty}\left\|A^{k _0}\right\| ^l \left\|\mathbf{C}\right\|K _A K _\omega=\frac{\left\|\mathbf{C}\right\|K _A K _\omega}{1-\left\|A^{k _0}\right\|}.
\end{equation*}
Finally, since each of the summands in the series is a continuous function then so is $f_{(\phi, \omega,F)}$.

\medskip

\noindent {\bf (ii)} The result that we just proved guarantees that if the differentials $D^if_{(\phi, \omega,F)}(m)$, $i \in  \left\{1, \ldots, r\right\}$, exist then they are given by the series $\sum _{j=0}^{\infty}A ^j \mathbf{C} D^i \left(\omega \circ \phi^{-j}\right) (m)$ that, using again the Weierstrass M-test and the hypotheses in the statement, will be now shown to uniformly converge to a continuous map. Indeed, using again the decomposition $A ^j=A^{l k _i+s}$ in terms of the element $k _i \in \mathbb{N}$ such that $\left\|A^{k _i}\right\| \left\|T ^i \phi ^{-k _i}\right\|_{\infty}<1   $ we can conclude that each summand of this series satisfies that
\begin{equation}
\label{bounds for derivatives}
\left\|A ^j \mathbf{C} D^i \left(\omega \circ \phi^{-j}\right)(m)\right\|\leq \left(\left\|A^{k _i}\right\| \left\|T \phi^{-k _i}\right\| _{\infty} \right)^l \left\|\mathbf{C}\right\|K _A^i K _{T^i \phi^{-1}} \left\|D^i \omega\right\|_{\infty},
\end{equation}
with $K _A^i:=\max \left\{1, \left\|A\right\|, \ldots, \left\|A^{k _i -1}\right\| \right\}$ and $K _{T^i \phi^{-1}}:= \max \left\{1,\left\|T^i \phi^{-1}\right\|_{\infty}, \ldots,  \left\|T^i \phi^{-1}\right\|_{\infty}^{k _i-1}\right\}$, which proves the desired convergence and that $D^if_{(\phi, \omega,F)}(m)=\sum _{j=0}^{\infty}A ^j \mathbf{C} D^i \left(\omega \circ \phi^{-j}\right) (m)$. Moreover,
\begin{equation}
\label{bounds for derivatives 2}
\left\|D^if_{(\phi, \omega,F)}(m)\right\|= \left\|\sum _{j=0}^{\infty}A ^j \mathbf{C} D^i \left(\omega \circ \phi^{-j}\right) (m)\right\|\leq \frac{\left\|\mathbf{C}\right\|K _A^i K _{T^i \phi^{-1}} \left\|D^i \omega\right\|_{\infty}}{1-\left\|A^{k _i}\right\| \left\|T^i \phi^{-k _i}\right\| _{\infty} }. 
\end{equation}

\medskip

\noindent {\bf (iii)} We start by noting that if the map \eqref{functor gs} is  continuous then the subsets $\Omega _i$ and $\Omega_e $ are indeed open because by Theorems 1.1 and 1.4 in \cite{Hirsch:book} the immersions and the embeddings in $C ^r(M, \mathbb{R}^N) $ are open and hence $\Omega _i$ and $\Omega_e $ are the preimages of those open sets by the continuous map $\Theta_{(\phi, F)} $. We establish now the continuity of $\Theta_{(\phi, F)} $ by showing that if the sequence $\left\{\omega _n\right\}_{n \in \mathbb{N}} $  in $C ^r(M, \mathbb{R}) $ converges to some element $\omega \in C ^r(M, \mathbb{R})  $ then so does $\left\{\Theta_{(\phi, F)} (\omega _n)\right\}_{n \in \mathbb{N}} \subset C ^r(M, \mathbb{R}^N) $ with respect to  $\Theta_{(\phi, F)} (\omega) \in C ^r(M, \mathbb{R}^N)  $. Indeed, if $\omega_n \longrightarrow \omega $ then, using the notation introduced in \eqref{bounds for derivatives 2}, we have that for a given $\epsilon> 0 $  and for $n$ sufficiently large 
\begin{equation*}
\frac{\left\|\mathbf{C}\right\|K _A^i K _{T^i \phi^{-1}}}{1-\left\|A^{k _i}\right\| \left\|T^i \phi^{-k _i}\right\| _{\infty} }  \left\|D^i \omega_n- D^i \omega\right\|_{\infty}< \epsilon/r.
\end{equation*}
Then,
\begin{multline*}
\left\|\Theta_{(\phi, F)} (\omega_n)- \Theta_{(\phi, F)} (\omega)\right\|_{C ^r(M, \mathbb{R}^N) }=\sum _{i=0} ^r \left\| D^if_{(\phi, \omega _n,F)}-  D^if_{(\phi, \omega ,F)} \right\| _{\infty}  \\
= \sum _{i=0} ^r
\left\|\sum _{j=0}^{\infty}A ^j \mathbf{C} D^i \left((\omega _n- \omega) \circ \phi^{-j}\right) (m)\right\|\leq
\sum _{i=0} ^r \frac{\left\|\mathbf{C}\right\|K _A^i K _{T^i \phi^{-1}}}{1-\left\|A^{k _i}\right\| \left\|T^i \phi^{-k _i}\right\| _{\infty} }  \left\|D^i \omega_n- D^i \omega\right\|_{\infty}<\frac{\epsilon}{r}+ \cdots+\frac{\epsilon}{r}= \epsilon,
\end{multline*}
as required. \quad $\blacksquare$

\section{Immersive generalized synchronizations}
\label{Immersive generalized synchronizations}

As we discussed in the introduction, the fact that the Takens delay map is an embedding under certain circumstances guarantees that the representation of the dynamical system associated to it can be used to learn the dynamics of its observations. In this section we take the first steps to show that similar results can be achieved by using the generalized synchronizations introduced in Proposition \ref{generalized synch with spectral radius}. More specifically, we shall spell out conditions on  linear state-space systems that guarantee that the resulting generalized synchronizations are immersions for generic scalar observations $\omega\in C^2(M, \mathbb{R})$. All along this section, the phase space manifold $M$ of the dynamical system $\phi \in {\rm Diff}^2(M)$ is compact and hence genericity in $C^2(M, \mathbb{R})$ is stated with respect to the topology associated to the extension to second-order differentiable functions of the Banach structure introduced in \eqref{definition c1d}. The next theorem is the main statement of this section.

\begin{theorem}
\label{Theorem immersion}
Let $\phi \in {\rm Diff}^2(M)$ be a dynamical system on a compact manifold $M$ of dimension $q$ that exhibits finitely many periodic orbits. Let $F: \mathbb{R}^N \times \mathbb{R} \longrightarrow \mathbb{R}^N  $ be a linear state map as in \eqref{linear system definition} with $N \geq 2q $ whose connectivity matrix satisfies that $\rho(A)<1 $ and such that for any observation map $\omega \in C ^2(M, \mathbb{R})$ the corresponding generalized synchronization  $f_{(\phi, \omega,F)} \in C ^2(M, \mathbb{R}^N) $ and, moreover, the map $\Theta_{(\phi, F)}: C ^2(M, \mathbb{R})  \longrightarrow  C ^2(M, \mathbb{R}^N)$ introduced in \eqref{functor gs} is continuous. Suppose also that the two following conditions hold:
\begin{description}
\item [(i)] For each periodic orbit $m$ of $\phi $ with period $n \in \mathbb{N} $, the derivative $T_m \phi^{-n} $ has $q$ distinct eigenvalues $\lambda_1, \lambda _2, \ldots, \lambda _q $. Let $\lambda _{{\rm max}} $ be the eigenvalue with the highest absolute value among the eigenvalues of all those linear maps and let $n_{ {\rm min}} $ be the smallest period. Suppose that $\rho(\lambda _{{\rm max}} A^{n_{ {\rm min}}})<1  $ and that for any periodic point $m$, the vectors
\begin{equation}
\label{condition on A C for later}
\left\{\left(\mathbb{I}- \lambda_j A ^n \right) ^{-1} \left(\mathbb{I} - A\right)^{-1} \left(\mathbb{I}- A ^n\right)\mathbf{C}\right\}_{j \in  \left\{1, \ldots, q\right\}}, \  \mbox{ with $\lambda_j$ eigenvalue of  $ T_m \phi^{-n}$} 
\end{equation}
form a linearly independent set.
\item [(ii)] The vectors $\left\{A ^j \mathbf{C}\right\}_{j \in \left\{0,1, \ldots, N-1 \right\}}$ form a linearly independent set.
\end{description}
Then, for generic $\omega \in C ^2(M, \mathbb{R}^N)$ the generalized synchronization $f_{(\phi, \omega,F)} \in C ^2(M, \mathbb{R}) $ is an immersion.
\end{theorem}

\paragraph{About the hypotheses of the theorem.}
All the hypotheses in this statement can be either easily guaranteed or, even better, they generically hold. More specifically, the condition on the linear state map $F$ to produce GS maps $f_{(\phi, \omega,F)} \in C ^2(M, \mathbb{R}^N) $ for any observation map $\omega \in C ^2(M, \mathbb{R})$ and the continuity of $\Theta_{(\phi, F)} $ can be enforced by using the second and third parts of Proposition \ref{generalized synch with spectral radius}. The condition on $\phi$ exhibiting finitely many periodic orbits holds generically due to the Kupka-Smale Theorem \cite{kupka1963contributiona, smale1963stable}.

As to the condition \eqref{condition on A C for later}, we shall see later on (see Proposition \ref{random_matrix_lemma}) that it holds almost surely in a very specific sense. Regarding the hypothesis in point {\bf (ii)}, this is a very important condition that amounts to reachability in a control theoretical sense  (see \cite{Kalman2010, sontag:book}). It has been shown in \cite{RC15} that if $A$ is diagonalizable then this condition holds if and only if all the eigenvalues in the spectrum $\sigma(A) $ of $A$ are distinct and in the linear decomposition $\mathbf{C}=\sum_{i=1}^N c _i \mathbf{v} _i$, with $\left\{\mathbf{v}_1, \ldots, \mathbf{v}_N\right\}$ a basis of eigenvectors of $A$, all the coefficients $c _i $, with $i \in \left\{1, \ldots, N\right\} $, are non-zero. This condition  can be equivalently reformulated by saying that the Krylov space \cite{krylov1931numerical} generated by $A$ and $\mathbf{C} $ has maximal dimension.

\paragraph{Relation with Takens' Theorem.} 
The system spelled out in the introduction that allows us to see Takens's  delay embedding
$\Phi_{(\phi, \omega)}$ as the GS  corresponding to a linear state map trivially satisfies the hypotheses of the theorem. Indeed, since in that case  $A$ is the lower shift matrix in dimension $2q+1  $ and $\mathbf{C}= (1,0, \ldots,0)^{\top} \in \mathbb{R}^{2q+1} $, the set in condition {\bf (ii)} coincides with the canonical basis in $ \mathbb{R}^N $ which is a trivially linearly independent set. Regarding the conditions in {\bf (i)}, as $A$ is nilpotent then all its eigenvalues are zero and hence the hypotheses are trivially satisfied.

Based on this observation, we can formulate a more general statement by saying that any linear system with nilpotent connectivity matrix $A$ that has an input mast $\mathbf{C} $ for which the vectors $\left\{A ^j \mathbf{C}\right\}_{j \in \left\{0,1, \ldots, N-1 \right\}}$ form a linearly independent set also satisfies the hypotheses of the theorem. Equivalently, with the terminology of the previous paragraph, we can rephrase this by writing that any reachable linear system with nilpotent connectivity matrix $A$ satisfies the hypotheses of the theorem.

\paragraph{System isomorphisms.} Given the linear state map introduced in \eqref{linear system definition} and a linear isomorphism of $\mathbb{R}^N $ with associated matrix $P \in \mathbb{M}_{N,N} $, consider the new map $\overline{F}(\mathbf{x}, z):=PAP ^{-1} \mathbf{x}+ P \mathbf{C}z  $. Let now $h: \mathbb{R}^M \longrightarrow \mathbb{R}^m  $ be a readout for the state map $F$. In this setup, it is easy to see that the state-space systems $(F,h)$  and $(\overline{F}, \overline{h}:=h \circ P ^{-1})$ are isomorphic in the sense that, in the presence of the echo state property, they determine identical input/output systems.

In view of this observation, it is important to emphasize that the hypotheses of Theorem \ref{Theorem immersion} are invariant under linear system isomorphisms. More explicitly, if we replace in the statement $A $  and $\mathbf{C}$ by $\overline{A}:=PAP ^{-1} $ and $\overline{\mathbf{C}}:=P \mathbf{C} $, respectively, then $\rho(A)= \rho(\overline{A}) $ and  the validity of the hypotheses {\bf (i)} and {\bf (ii)} is not altered. Indeed, regarding {\bf (i)}, it suffices to notice that
\begin{multline*}
\left(\mathbb{I}- \lambda_j \overline{A} ^n \right) ^{-1} \left(\mathbb{I} - \overline{A}\right)^{-1} \left(\mathbb{I}- \overline{A} ^n\right)\overline{\mathbf{C}}=\left(\sum _{i=0}^{\infty}\lambda_j ^i \left(\overline{A} ^n\right)^i\right)\left(\sum _{i=0}^{\infty}\lambda_j ^i \overline{A} ^i\right)\left(\mathbb{I}- \overline{A} ^n\right)\overline{\mathbf{C}}\\
=P\left(\sum _{i=0}^{\infty}\lambda_j ^i \left({A} ^n\right)^i\right)P ^{-1}P\left(\sum _{i=0}^{\infty}\lambda_j ^i \overline{A} ^i\right)P ^{-1}\left(PP ^{-1}- P{A} ^nP ^{-1}\right)P{\mathbf{C}}\\=P \left(\left(\mathbb{I}- \lambda_j A ^n \right) ^{-1} \left(\mathbb{I} - A\right)^{-1} \left(\mathbb{I}- A ^n\right)\mathbf{C}\right).
\end{multline*}
As to {\bf (ii)}, notice that $\left\{\overline{A} ^j \overline{\mathbf{C}}\right\}_{j \in \left\{0,1, \ldots, N-1 \right\}}=\left\{PA ^j \mathbf{C}\right\}_{j \in \left\{0,1, \ldots, N-1 \right\}}$. Since $P$ is an invertible matrix, in both cases the linear independence is preserved.

Another observation that is worth pointing out is that the class of linear systems for which Theorem \ref{Theorem immersion} hold is strictly larger than the one determined (up to linear isomorphisms) by Takens' Theorem. As it was mentioned in the previous paragraph, Takens' result is associated to a linear system with nilpotent connectivity matrix $A$ (whose eigenvalues are hence all zero). It is easy to see that when the entries of $\mathbf{C}$ are all non-zero then one can always find a non-singular diagonal matrix $A$ for which the hypotheses of Theorem \ref{Theorem immersion} hold. Such system is not in the same isomorphism class as Takens' system.

\paragraph{Proof of the Theorem.} We proceed in two steps. In the first one we show that $f_{(\phi, \omega,F)} \in C ^2(M, \mathbb{R}) $ is an immersion at periodic points and in the second one we take care of the remaining points. We emphasize that equilibria can be seen as periodic points with period $1$.

\medskip

\noindent {\bf Step 1. Immersion at periodic points.} We start this part with two preparatory lemmas.

\begin{lemma}
Consider a connectivity matrix that satisfies the conditions $\rho(A)<1 $ and also that $\rho(\lambda _{{\rm max}} A^{n_{ {\rm min}}})<1$ as in part {\bf (i) } of the statement of the theorem. Then, for any periodic point $m$ with period $n$ and any  eigenvalue $\lambda_j$ of  $ T_m \phi^{-n}$, we have that $\rho(\lambda _j A^{n})<1$ and 
\begin{equation}
\label{funny inverse for per}
(\mathbb{I} - \lambda_j A^n)^{-1}=\sum_{k=0} ^{\infty} \lambda _j^kA^{nk}.
\end{equation}
\end{lemma}

\noindent\textbf{Proof.\ \ } Firstly, recall the general fact already used in the proof of Proposition \ref{generalized synch with spectral radius} (see also Proposition 4.2 in \cite{RC16}) that for any square matrix  $B$ such that $\rho(B)< 1  $ then $\left(\mathbb{I}-B\right) ^{-1}=\sum_{j=0}^{\infty} B ^j $. Let now $m$ be a periodic point with period $n$ and let $\lambda_j$ be an  eigenvalue  of  $ T_m \phi^{-n}$. 
 This implies that in order for \eqref{funny inverse for per} to hold we just need to show that $\rho(\lambda _j A^{n})<1$. This is indeed true since  any element in the spectrum of $\lambda _j A^{n} $ can be written as $\lambda _j \mu _k ^n $ with $\mu_k \in \mathbb{C}$ an eigenvalue of $A$. Moreover, let $c < 1$ such that $\left| \lambda _j \right|=c \left| \lambda _{{\rm max}}  \right| $. Then
\begin{equation*}
\left| \lambda _j \mu _k ^n \right|=c\left|\lambda _{{\rm max}} \mu _k ^{n_{{\rm min}}}  \right|\left|  \mu _k ^{n-n_{{\rm min}}} \right|< 1, 
\end{equation*} 
as required. Notice that in the last inequality we used that $\rho(\lambda _{{\rm max}} A^{n_{ {\rm min}}})<1$ and that $\rho(A)<1 $. $\blacktriangledown $.

\begin{lemma}
\label{condition for linind}
In the hypotheses of the statement of the theorem, let $m \in M $ be a periodic point of $\phi \in {\rm Diff}^2(M)$ with period $n \in \mathbb{N} $. Let $\left\{\mathbf{v} _1, \ldots, \mathbf{v} _p\right\} $ be a basis of eigenvectors associated to the distinct eigenvalues $\lambda_1, \lambda _2, \ldots, \lambda _q $.  Suppose that the set
\begin{equation}
\label{linind condition}
\left\{(\mathbb{I} - \lambda_j A^n)^{-1} \sum^{n-1}_{k = 0} A^k \mathbf{C} D (\omega \circ  \phi^{-k})(m) \mathbf{v}_j\right\}_{j \in \left\{1, \ldots,q\right\}}
\end{equation}
is linearly independent.
Then $f_{(\phi, \omega,F)}$ is an immersion at the periodic point $m$ for generic $\omega \in C ^2(M, \mathbb{R}^N)$.
\end{lemma}

\noindent\textbf{Proof.\ \ } Since the eigenvalues $\lambda_j$ are distinct and the eigenvectors $\mathbf{v}_j$ are hence linearly independent, it therefore suffices to show that the set $\left\{Df_{(\phi, \omega,F)}(m) \mathbf{v}_j\right\}_{j \in \left\{1, \ldots,q\right\}}$ is linearly independent to conclude that $Df_{(\phi, \omega,F)}(m)$ is injective. Then, by the expression \eqref{higher order derivatives of GS}:
\begin{align*}
Df_{(\phi, \omega,F)}(m)\mathbf{v} _j &= \sum^{\infty}_{l = 0} A^l \mathbf{C} D (\omega \circ \phi^{-l})(m) \mathbf{v} _j= \sum^{\infty}_{l = 0} \sum_{k = 0}^{n-1} A^{ln+k} \mathbf{C} D(\omega \circ \phi^{-(ln+k)})(m)\mathbf{v} _j \\
&= \sum^{\infty}_{l = 0} \sum_{k = 0}^{n-1} A^{ln+k} \mathbf{C} D(\omega \circ \phi^{-k})(m) [D\phi^{-n}(m)]^l \mathbf{v} _j = \sum^{\infty}_{l = 0} \sum_{k = 0}^{n-1} A^{ln+k} \mathbf{C} D(\omega \circ  \phi^{-k})(m) \lambda^l_j \mathbf{v}_j \\
&= \sum^{\infty}_{l = 0} (\lambda_jA^n)^l \sum_{k = 0}^{n-1} A^k \mathbf{C} D(\omega \circ  \phi^{-k})(m)  \mathbf{v}_j = 
(\mathbb{I} - \lambda_j A^n)^{-1} \sum^{n-1}_{k = 0} A^k \mathbf{C} D (\omega \circ  \phi^{-k})(m) \mathbf{v}_j,
\end{align*}
which proves the statement. $\blacktriangledown $.

\medskip

We now use this result to show that, for generic $\omega \in C ^2(M, \mathbb{R}^N)$, the generalized synchronization $f_{(\phi, \omega,F)} \in C ^2(M, \mathbb{R}) $ is an immersion at the periodic points of $\phi$. Let $m _1, \ldots, m _P \in M$ be the distinct periodic points of $\phi$, each of which have periods $n _1, \ldots, n _P \in \mathbb{N} $, respectively (the equilibria of $\phi$ are on this list with periods equal to one). The term {\it distinct} means that none of those points are in the orbits of the others. We now choose $P$ disjoint open neighborhoods $B _i $ that contain each of the distinct periodic points $m _i $. Since there is a finite number of periodic points, the open sets $B _i $ can be chosen small enough so that, additionally, all the open sets
\begin{equation*}
\phi^{-t}(B _i) \quad \mbox{for all $t \in \left\{0, \ldots, n _i\right\} $ and $i \in \left\{1, \ldots, P\right\}$}
\end{equation*}
are disjoint.

Now, given any of the distinct periodic points $ m _i \in M $ on the list, we show that  $f_{(\phi, \omega,F)} $ for generic $\omega \in C ^2(M, \mathbb{R}^N)$, that is, the set of observation maps $\omega$ for which $f_{(\phi, \omega,F)} $ is an immersion at $m _i $ is open and dense in $C ^2(M, \mathbb{R}^N)$. The openness is a consequence of the hypothesis on the continuity of the map $\Theta_{(\phi, F)} $ and of an argument identical to the beginning of the proof of part {\bf (iii)} of Proposition \ref{generalized synch with spectral radius}. Regarding the density, we show that  if $f_{(\phi, \omega,F)} $ is not an immersion at $m _i  $, then there is a perturbation $\omega' $ of $\omega $ in $C^2(M, \mathbb{R})$ for which $f_{(\phi, \omega',F)} $ is an immersion at $m _i $. Indeed, set 
\begin{equation}
\label{perturbed observation map 1}
    \omega' = \omega + \sum_{l = 0}^{n _i-1}   \psi_l ^i
\end{equation}
where  $\psi_l ^i \in C^{\infty}(M,\mathbb{R})$ are bump functions whose  supports are contained in  $\phi^{-l}(B_i)$ and, additionally, are chosen to satisfy
\begin{align*}
D(\psi_l^i \circ \phi^{-l})(m _i) = \varepsilon \mathbf{v}^{\top}, \quad \mbox{$l \in \left\{0, \ldots, n _i-1\right\}$,}
\end{align*}
for some small constant $\varepsilon> 0$ and $\mathbf{v} \in \mathbb{R} ^p $  the unique vector that solves the linear system
\begin{equation}
\label{linear system for later}
\left(
\begin{array}{c}
\mathbf{v}_1^{\top}\\
\vdots\\
\mathbf{v}_p^{\top}
\end{array}
\right) \mathbf{v}=
\left(
\begin{array}{c}
1\\
\vdots\\
1
\end{array}
\right),
\end{equation}
with $\left\{\mathbf{v} _1, \ldots, \mathbf{v} _p\right\} $ a basis of eigenvectors of $T_{m _i}\phi^{-n _i} $. Note that by construction and for any $l \in \left\{0, \ldots, n _i-1\right\}$,
\begin{equation}
\label{another inter 1}
D(\omega'\circ \phi^{-l})(m _i) = D(\omega\circ \phi^{-l})(m _i) + D(\psi_l^i \circ  \phi^{-l}) = D(\omega\circ \phi^{-l})(m _i) + \varepsilon \mathbf{v}^{\top}.
\end{equation}
We now consider the vectors \eqref{linind condition} in Lemma \ref{condition for linind} with respect to the perturbed observation map in \eqref{perturbed observation map 1}. Indeed, by \eqref{another inter 1} and the way in which the vector $\mathbf{v}  $ has been chosen in \eqref{linear system for later}:
\begin{multline*}
(\mathbb{I} - \lambda_j A^{n_i})^{-1} \sum^{n_i-1}_{k=0} A^k \mathbf{C} D(\omega' \circ \phi^{-k})(m _i)\mathbf{v}_j \\
= (\mathbb{I} - \lambda_j A^{n_i})^{-1} \sum^{n_i-1}_{k=0} A^k \mathbf{C} D(\omega \circ  \phi^{-k})(m _i)\mathbf{v}_j + \varepsilon(\mathbb{I} - \lambda_j A^{n_i})^{-1} \sum^{n_i-1}_{k=0} A^k \mathbf{C}  \mathbf{v}^{\top} \mathbf{v}_j \\
= (\mathbb{I} - \lambda_j A^{n_i})^{-1} \sum^{n_i-1}_{k=0} A^k \mathbf{C} D(\omega \circ  \phi^{-k})(m _i)\mathbf{v}_j 
+ \varepsilon(\mathbb{I} - \lambda_j A^{n_i})^{-1} (\mathbb{I} - A)^{-1}(I - A^{n_i})\mathbf{C}.
\end{multline*}
Given that when we vary $j \in \left\{1, \ldots, p \right\}$ in the previous expression the vectors in the second summand form by hypothesis a linearly independent set, we can use Lemma \ref{invertible perturbation} to choose $\varepsilon> 0 $ so that the family $\left\{(\mathbb{I} - \lambda_j A^{n_i})^{-1} \sum^{n_i-1}_{k=0} A^k \mathbf{C} D(\omega' \circ \phi^{-k})(m _i)\mathbf{v}_j\right\}_{j \in \left\{1, \ldots, p \right\}}$ forms a linearly independent set and, at the same time, $\omega'$ is as close to $\omega $ in $C^2(M, \mathbb{R})$ as desired. This shows by Lemma \ref{condition for linind} that $f_{(\phi, \omega',F)} $ is an immersion at $m _i $. 

The choice of the open sets $B _i $ implies that we can keep perturbing $\omega$ in order to make $f_{(\phi, \omega',F)} $ immersive at the other periodic points without spoiling that condition for the previous ones. This shows in particular that a perturbation of the type
\begin{equation}
\label{perturbed observation map final}
    \omega' = \omega + \sum_{i=1}^P\sum_{l = 0}^{n _i-1}   \psi_l ^i
\end{equation}
can be constructed so that $f_{(\phi, \omega',F)} $ is immersive at all the periodic points of $\phi $, as required.

\medskip

\noindent {\bf Step 2. Immersion at the remaining points.} Having just proved that for generic $\omega \in C ^2(M, \mathbb{R}^N)$ the generalized synchronization $f_{(\phi, \omega,F)} \in C ^2(M, \mathbb{R}) $ is an immersion at the periodic points of $\phi$, the Immersion Theorem (see \cite[Theorem 3.5.7]{mta}) guarantees that the same holds for the open set  formed by the union of certain open neighborhoods around those points. 
Let $\mathcal{M} \subset M$ be the compact subset of $M$ obtained by removing that immersed open set. Our goal is now to show that $f_{(\phi, \omega,F)} \in C ^2(M, \mathbb{R}) $ is also an immersion at $\mathcal{M}$ for generic $\omega \in C ^2(M, \mathbb{R}^N)$. 

Recall first that the hypotheses that we imposed on $M$ in the beginning of Section \ref{Definitions and preliminary discussion} imply that it can be endowed with a Riemannian metric which makes it into a complete metric space by the Hopf and Rinow Theorem (see \cite[Theorem 7.7]{boothby2003introduction}). This implies in turn that the compact subset ${\cal M}\subset M $ is also a complete metric space which allows us to define open balls $B _r(m) $ of radius $r > 0$  around each point $m \in {\cal M} $. Using this notation, in the next paragraphs we show that for any $\omega \in C ^2(M, \mathbb{R}^N)$ and $m  \in {\cal M} $ we can find a $n (m) \in \mathbb{N}  $  and a perturbation $\omega' \in C ^2(M, \mathbb{R}^N)$ as close to $\omega$ as desired such that the restriction of  $f_{(\phi, \omega',F)} $ to $B _{2^{-n(m)}}(m) $ is an immersion.

Indeed, take an arbitrary $m \in \mathcal{M}$ and define a collection of balls $B_{2^{-n}}(m)$ centered at $m$ with radius $2^{-n}$,  $n \in \mathbb{N}$. For a fixed $n$ consider the infinite trajectory $\phi^{-t}(B_{2^{-n}}(m))$, $t  \in \mathbb{N}$.  Choose now $n_1 \in \mathbb{N}$ large enough so that, for any $n > n_1$ the balls $\phi^{-t}(B_{2^{-n}}(m))$   are disjoint for $t = 0, \ldots , N-1$ and $B_{2^{-n}}(m) \subset U$ where $(U, h) $ is an admissible chart of $M$. Given that $\phi \in {\rm Diff}^2(M)$, we note that the family $(U _t, h _t)$, $t \in \mathbb{N} $,  defined by $U _t= \phi^{-t}(U) $ and  $h _t :=h \circ \phi ^t  $ is made of admissible charts and that $\phi^{-t}(B_{2^{-n}}(m))\subset U_t$, for all $n> n _1 $. Let $T(n)$ denote the largest integer such that $\phi^{-t}(B_{2^{-n}}(m))$ are disjoint  for $t = 0, \ldots, T (n)-1$.

Now, for each $n > n _1$ and $t = 0, \ldots, N-1$ we define functions $\psi_{tn} \in C^{\infty}(M,\mathbb{R})$ that have their support included in $\phi^{-t}(B_{2^{-n}}(m))$ and satisfy
\begin{equation}
\label{condition on the derivative}
\frac{\partial (\psi_{tn}h_t^{-1})}{\partial u_j} = 1
\end{equation}
on $h_t\left(\phi^{-t}(B_{2^{-(n+1)}}(m))\right) =h\left( B_{2^{-(n+1)}}(m)\right)$. We impose further that $\psi_{tn} = \psi_{t(n+1)}$ on $\phi^{-t}(B_{2^{-(n+2)}}(m))$ for all $n> n _1$, and that there is some $\kappa > 0$ independent of $n$ and $t$ such that $\lVert \psi_{tn} h_t^{-1} \rVert_{C^1} \leq \kappa$.
These functions can be constructed by setting
\begin{align*}
    \psi_{tn}(m) = \lambda_{tn}(m) \sum_{j = 1}^{q}\xi_j(m)
\end{align*}
where $\xi_j$ is the $j$-th coordinate map for the chart $h _t$ and $\lambda_{tn} \in C^{\infty}(M,\mathbb{R})$ are bump functions that have support included in $\phi^{-t}(B_{2^{-n}}(m))$ and satisfy $\lambda_{tn}|_{\phi^{-t}(B_{2^{-(n+1)}}(m))} = 1$. Define now the perturbation $\omega _n $ of $\omega$ by
\begin{align}
\label{definition of omega perturb}
    \omega_n = \omega + \sum_{t = 0}^{N-1} \varepsilon_t \psi_{tn},
\end{align}
where $\varepsilon _t  $  are the components of a vector $\boldsymbol{\varepsilon} \in \mathbb{R}^N  $ with positive entries. By construction, for any $m' \in B_{2^{-n}}(m)$ and $t = 0 , \ldots,  N-1$, we have  that
\begin{align*}
    \omega_n \phi^{-t}(m') = \omega \phi^{-t}(m') + \varepsilon_t \psi_{tn}(m') 
\end{align*}
and moreover by \eqref{condition on the derivative} and for any $m' \in B_{2^{-(n+1)}}(m)$:
\begin{equation}
\label{component of delay}
    \frac{\partial (\omega_n \phi^{-t} h^{-1})}{\partial u_j}(h(m')) = 
    \frac{\partial (\omega \phi^{-t} h^{-1})}{\partial u_j}(h(m')) + \varepsilon_t.
\end{equation}
Let $\Phi: M \to \mathbb{R}^N$ be a backwards version of the Takens delay map introduced in \eqref{delay embedding expression}, that is, 
\begin{align*}
    \Phi(m) := \left( \omega(m) , \omega\circ  \phi^{-1}(m) , \ldots, \omega\circ\phi^{-(N-1)}(m) \right)^{\top},
\end{align*}
and let $\Phi_n : M \to \mathbb{R}^N$ be its perturbed version defined by
\begin{align*}
    \Phi_n(m) := \left( \omega_n(m) , \omega_n\circ \phi^{-1}(m) , \ldots, \omega_n\circ\phi^{-(N-1)}(m) \right)^{\top}.
\end{align*}
Using these objects, we can rewrite \eqref{component of delay} in vector form as
\begin{equation}
\label{derivative in vector form bis}
\frac{\partial (\Phi_n h^{-1})}{\partial u_j}(h(m')) = \frac{\partial (\Phi h^{-1})}{\partial u_j} (h(m'))+ \boldsymbol{\varepsilon},
\end{equation}
for any $m' \in B_{2^{-(n+1)}}(m)$ and where $\boldsymbol{\varepsilon} \in \mathbb{R}^N$. Next, for any $t = N , \ldots, T(n) - 1$ notice that
\begin{equation}
\label{derivative in vector form bis1}
    \omega_n \phi^{-t}(m) = \omega \phi^{-t}(m).
\end{equation}
Finally, if $t \geq T(n)$ then 
\begin{equation}
\label{derivative in vector form bis2}
    \omega_n \phi^{-t}(m) = \omega \phi^{-t}(m) + \sum^{N-1}_{\tau = 0} \varepsilon_\tau \psi_{\tau n}\phi^{-t}(m).
\end{equation} 
We now consider the perturbed generalized synchronization $f_{(\phi, \omega_n,F)} :M \longrightarrow\mathbb{R} ^N$  given by
\begin{align*}
    f_{(\phi, \omega_n,F)} &= \sum_{t = 0}^{\infty} A^t \mathbf{C} \omega_n \phi^{-t} = \sum_{t = 0}^{N-1} A^t \mathbf{C} \omega_n \phi^{-t} + \sum_{t = N}^{\infty} A^t \mathbf{C} \omega_n \phi^{-t} = Q \Phi_n + \sum_{t = N}^{\infty} A^t \mathbf{C} \omega_n\phi^{-t}
\end{align*}
where $Q$ is the $N \times N$ real matrix with $(t+1)$-th column $A^t \mathbf{C}$. Now we take the partial derivatives with respect to $u_j$ at points in $h(B_{2^{-(n+1)}}(m))$ and observe that by \eqref{derivative in vector form bis}, \eqref{derivative in vector form bis1}, and \eqref{derivative in vector form bis2}:	
\begin{align}
    \frac{\partial (f_{(\phi, \omega_n,F)} h^{-1})}{\partial u_j} &= Q \frac{\partial(\Phi_n h^{-1})}{\partial u_j} + \sum_{t = N}^{\infty} A^t \mathbf{C} \frac{\partial (\omega_n\phi^{-t}h^{-1} )}{\partial u_j} =Q \frac{\partial(\Phi h^{-1})}{\partial u_j} + Q \boldsymbol{\varepsilon} + \sum^{\infty}_{t = N} A^t \mathbf{C} \frac{\partial (\omega_n \phi^{-t}h^{-1})}{\partial u_j}\notag \\
    &=
    Q \frac{\partial(\Phi h^{-1})}{\partial u_j} + Q \boldsymbol{\varepsilon} + \sum^{T(n)-1}_{t = N} A^t \mathbf{C} \frac{\partial (\omega_n \phi^{-t} h^{-1})}{\partial u_j} + 
    \sum^{\infty}_{t = T(n)} A^t \mathbf{C} \frac{\partial (\omega_n \phi^{-t} h^{-1})}{\partial u_j}\notag \\
    &= Q \frac{\partial(\Phi h^{-1})}{\partial u_j} + Q \boldsymbol{\varepsilon} + \sum^{T(n)-1}_{t = N} A^t \mathbf{C} \frac{\partial (\omega \phi^{-t} h^{-1})}{\partial u_j} + 
    \sum^{\infty}_{t = T(n)} A^t \mathbf{C} \frac{\partial (\omega_n \phi^{-t} h^{-1})}{\partial u_j}\notag \\
    &= Q \frac{\partial(\Phi h^{-1})}{\partial u_j} + Q \boldsymbol{\varepsilon} + \sum^{T(n)-1}_{t = N} A^t \mathbf{C} \frac{\partial (\omega \phi^{-t}h^{-1})}{\partial u_j} + 
    \sum^{\infty}_{t = T(n)} A^t \mathbf{C} \frac{\partial (\omega \phi^{-t}h^{-1})}{\partial u_j}\notag \\
    &+ \sum^{\infty}_{t = T(n)} A^t \mathbf{C} \bigg( \sum^{N-1}_{\tau = 0} \varepsilon_{\tau} \frac{\partial( \psi_{\tau n} \phi^{-t} h^{-1})}{\partial u_j} \bigg)\notag \\
    &= Q \frac{\partial(\Phi h^{-1})}{\partial u_j} + Q \boldsymbol{\varepsilon} + \sum^{\infty}_{t = N} A^t \mathbf{C} \frac{\partial( \omega \phi^{-t}h^{-1})}{\partial u_j} + \sum^{\infty}_{t = T(n)} A^t \mathbf{C} \bigg( \sum^{N-1}_{\tau = 0} \varepsilon_{\tau} \frac{\partial (\psi_{\tau n} \phi^{-t} h^{-1})}{\partial u_j} \bigg)\notag \\
    &= \sum_{t = 0}^{N-1}A^t \mathbf{C} \frac{\partial (\omega \phi^{-t}h^{-1})}{\partial u_j} + Q \boldsymbol{\varepsilon} + \sum^{\infty}_{t = N} A^t \mathbf{C} \frac{\partial (\omega \phi^{-t}h^{-1})}{\partial u_j} + \sum^{\infty}_{t = T(n)} A^t \mathbf{C} \bigg( \sum^{N-1}_{\tau = 0} \varepsilon_{\tau} \frac{\partial (\psi_{\tau n} \phi^{-t} h^{-1})}{\partial u_j} \bigg)\notag \\
    &= \sum_{t = 0}^{\infty}A^t \mathbf{C} \frac{\partial (\omega \phi^{-t}h^{-1})}{\partial u_j} + Q \boldsymbol{\varepsilon} + \sum^{\infty}_{t = T(n)} A^t \mathbf{C} \bigg( \sum^{N-1}_{\tau = 0} \varepsilon_{\tau} \frac{\partial (\psi_{\tau n} \phi^{-t} h^{-1})}{\partial u_j} \bigg)\notag \\
    &= \frac{\partial(f_{(\phi, \omega,F)} h^{-1})}{\partial u_j} + Q \boldsymbol{\varepsilon} + \sum^{\infty}_{t = T(n)} A^t \mathbf{C} \bigg( \sum^{N-1}_{\tau = 0} \varepsilon_{\tau} \frac{\partial (\psi_{\tau n} \phi^{-t} h^{-1})}{\partial u_j} \bigg)
\label{important identity for derivatives}.
\end{align}
In order to prove that  $f_{(\phi, \omega,F)}$ is an immersion  at the  points in $h(B_{2^{-(n+1)}}(m))$ for a generic observation $\omega$, we shall find an arbitrarily small vector $\boldsymbol{\varepsilon}$ for which the vectors corresponding to the $\boldsymbol{\varepsilon}$-perturbed observation $\omega _n$ 
\begin{align*}
    \left\{\frac{\partial \left(f_{(\phi, \omega_n,F)}  h^{-1}\right)}{\partial u_j} \right\}_{j \in  \left\{1, \ldots, q\right\}}
\end{align*}
are a linearly independent family. We will proceed inductively by showing that if we assume for some $s$ satisfying $1 \leq s < q$ that the vectors
\begin{align}
\label{proceeding_vectors}
    \left\{\frac{\partial \left(f_{(\phi, \omega,F)}  h^{-1}\right)}{\partial u_j} \right\}_{j \in  \left\{1, \ldots, s\right\}}
\end{align}
are linearly independent, then we can choose an arbitrarily small vector $\boldsymbol{\varepsilon}$ such that the family corresponding to the perturbed observation $\omega_n$ defined in  \eqref{definition of omega perturb} satisfies that 
\begin{align*}
    \left\{\frac{\partial \left(f_{(\phi, \omega_n,F)}  h^{-1}\right)}{\partial u_j} \right\}_{j \in  \left\{1, \ldots, s+1\right\}}
\end{align*}
is a linearly independent family.  To this end, we define the map $\Psi : \mathbb{R}^s \times h(U) \to \mathbb{R}^{N}$ as
\begin{align*}
    \Psi(\boldsymbol{\alpha} , {\bf u}) = \sum^s_{j = 1}\alpha_j \frac{\partial \left(f_{(\phi, \omega,F)}  h^{-1}\right)}{\partial u_j} -  \frac{\partial \left(f_{(\phi, \omega,F)}  h^{-1}\right)}{\partial u_{s+1}}.
\end{align*}
The hypothesis on the statement of the theorem about $f_{(\phi, \omega,F)} \in C ^2(M, \mathbb{R}^N) $ for any observation map $\omega \in C ^2(M, \mathbb{R})$ implies that $\Psi $ is of class $C^1$ and maps a manifold  of dimension  $s + q$ to a manifold of dimension $N$. 
Since by hypothesis $s + q< 2q\leq N$,  then the set  $\mathbb{R}^N \setminus   \Psi(\mathbb{R}^s \times h(U))$ is dense in $\mathbb{R}^N  $ (see \cite[Chapter 3, Proposition 1.2]{Hirsch:book}). This implies that we can choose an arbitrarily small vector 
$ \boldsymbol{\delta} \in \bigg(\mathbb{R}^N \setminus \Psi(\mathbb{R}^s \times h (U))\bigg)$ such that  if we set $\boldsymbol{\varepsilon} := Q^{-1}\boldsymbol{\delta} $ then we have that the vector
\begin{align*}
    \frac{\partial(f_{(\phi, \omega,F)} h^{-1})}{\partial u_{s+1}} + Q \boldsymbol{\varepsilon}
\end{align*}
is independent of the vectors in \eqref{proceeding_vectors} when evaluated at any point in $h(U) $. Since the linear independence is stable under small perturbations we can choose $\boldsymbol{\varepsilon} $ small enough so that it is actually the family
\begin{align*}
    \left\{\frac{\partial \left(f_{(\phi, \omega,F)}  h^{-1}\right)}{\partial u_j}+ Q \boldsymbol{\varepsilon} \right\}_{j \in  \left\{1, \ldots, s+1\right\}}
\end{align*}
that is linearly independent. Now, in view of the identity \eqref{important identity for derivatives} we note that the value $n \in \mathbb{N}$ can be chosen large enough so that the residual terms
\begin{align*}
    \sum^{\infty}_{t = T(n)} A^t \mathbf{C} \bigg( \sum^{N-1}_{\tau = 0} \varepsilon_{\tau} \frac{\partial (\psi_{\tau n} \phi^{-t} h^{-1})}{\partial u_j} \bigg), \quad \mbox{$j \in  \left\{1, \ldots, s+1\right\}$,}
\end{align*} 
are small enough so that the family 
\begin{align*}
    \left\{\frac{\partial \left(f_{(\phi, \omega,F)}  h^{-1}\right)}{\partial u_j}+ Q \boldsymbol{\varepsilon}+  \sum^{\infty}_{t = T(n)} A^t \mathbf{C} \bigg( \sum^{N-1}_{\tau = 0} \varepsilon_{\tau} \frac{\partial (\psi_{\tau n} \phi^{-t} h^{-1})}{\partial u_j} \bigg)\right\}_{j \in  \left\{1, \ldots, s+1\right\}}=
  \left\{\frac{\partial \left(f_{(\phi, \omega_n,F)}  h^{-1}\right)}{\partial u_j}\right\}_{j \in  \left\{1, \ldots, s+1\right\}}
\end{align*}
is linearly independent, as required. Notice that this equality is a consequence of \eqref{important identity for derivatives}. The possibility to shrink the residual term comes from the convergence of the series
\begin{align*}
    \sum^{\infty}_{t = 0} A^t \mathbf{C} \bigg( \sum^{N-1}_{\tau = 0} \varepsilon_{\tau} \frac{\partial (\psi_{\tau n} \phi^{-t} h^{-1})}{\partial u_j} \bigg), \quad \mbox{$j \in  \left\{1, \ldots, s+1\right\}$,}
\end{align*} 
which is guaranteed by the hypothesis on the differentiability of $f_{(\phi, \omega,F)} $ for any observation map $\omega \in C ^2(M, \mathbb{R})$ and the expression \eqref{higher order derivatives of GS}. In this case the bump functions play the role of the observations for which we assumed the existence of a uniform bound $\kappa$ over $\tau$ and $n$ such that $\lVert \psi_{\tau n} \rVert_{C^1} < \kappa$.

If we recursively apply this procedure, we can conclude the existence of a small perturbation $\omega_n $ of $\omega $ obtained as a sequence of perturbations of the type \eqref{definition of omega perturb} for which the family 
\begin{equation*}
  \left\{\frac{\partial \left(f_{(\phi, \omega_n,F)}  h^{-1}\right)}{\partial u_j}\right\}_{j \in  \left\{1, \ldots, q\right\}}
\end{equation*}
is linearly independent when evaluated at $m \in {\cal M} $, which proves that 
$f_{(\phi, \omega_n,F)} \in C ^2(M, \mathbb{R} ^N ) $ is an immersion at $m \in {\cal M}$.

Finally, observe that we just showed that for any $m \in \mathcal{M}$, there exists an $n(m) \in \mathbb{N}$ such that the restriction of the perturbation $f_{(\phi, \omega_{n(m)},F)} $ to $B_{2^{-n(m)}}(m)$ is an immersion. We note that the union 
\begin{align*}
    \bigcup_{m \in \mathcal{M}} B_{2^{-n}}(m)
\end{align*}
is clearly an open cover of $\mathcal{M}$. Since $\mathcal{M}$ is compact, it admits a finite subcover. The finite subcover comprises sets for which, one at a time, we can construct an immersion using the procedure described earlier in this proof. For each set, we ensure that the perturbation is sufficiently small not to spoil the immersion on any other set. 

This argument completes the proof of the  immersion of the GS at the points of $\mathcal{M}$ and therefore, together with the Step 1, shows that there exists a small perturbation of $f_{(\phi, \omega_n,F)} \in C ^2(M, \mathbb{R} ^N ) $ of $f_{(\phi, \omega,F)}$ that is an immersion at all the points in $M$, as required. \quad $\blacksquare$

\section{Linear reservoir embeddings}
\label{Linear reservoir embeddings}

We continue in this section by showing two important facts. Firstly, we prove that without additional hypotheses, the globally immersive generalized synchronizations whose existence we proved in Theorem \ref{Theorem immersion} are injective and hence are necessarily embeddings due to the compactness of $M$ (see \cite{Hirsch:book}). As we already pointed out in the introduction this is very important in relation to the {\it learnability question}, that is, at the time of using the embedded state representation of the dynamical system to learn from data the dynamics of its observations. The second fact is related with the reservoir computing phenomenon as we show that randomly generated linear systems yield  synchronization maps $f_{(\phi, \omega_n,F)} \in C ^2(M, \mathbb{R} ^N ) $ that are almost surely embeddings and are hence amenable to learnability from data.

\begin{theorem}
\label{the immersions are embeddings}
Assume that the hypotheses of Theorem \ref{Theorem immersion} hold true and that, additionally, $N>\max \left\{2q, \ell\right\}$ with $\ell \in \mathbb{N} $ the lowest common multiple of all the periods of the finite periodic points of $\phi$. Then, for generic $\omega \in C ^2(M, \mathbb{R}^N)$, the generalized synchronization $f_{(\phi, \omega,F)} \in C ^2(M, \mathbb{R}^N) $ is an embedding.
\end{theorem}

\noindent\textbf{Proof.\ \ } As in the previous theorem, we proceed in two steps.

\medskip

\noindent {\bf Step 1. Injectivity around the periodic set.} We start by showing that the observations corresponding to the globally immersive generalized synchronizations whose existence we proved in Theorem \ref{Theorem immersion} can be slightly perturbed in $ C ^2(M, \mathbb{R})$ so that the resulting GS is injective in an open subset $V _P  $ that includes all the periodic points of $\phi$. We start this part of the discussion with a preparatory lemma.

\begin{lemma}
\label{injectivity_lemma}
In the hypotheses of the theorem, let $m _1, \ldots, m _P \in M$ be the distinct periodic points of $\phi$, each of which have periods $n _1, \ldots, n _P \in \mathbb{N} $, respectively. Let $\ell \in \mathbb{N} $ be the lowest common multiple of all the periods and denote by $M _P $ the set of all periodic points of $\phi$ (that is, the set that comprises $ \left\{m _1, \ldots, m _P\right\} $ and all the corresponding orbits). Then, the restriction $\left.f_{(\phi, \omega,F)}\right|_{M _P} $ of a generalized synchronization $f_{(\phi, \omega,F)} \in C ^2(M, \mathbb{R}^N) $ to $M _P $  is injective if and only if the map $g_{(\phi, \omega,F)} : M _P \longrightarrow \mathbb{R}^N$ defined by
    \begin{align*}
        g_{(\phi, \omega,F)} = \sum_{k = 0}^{\ell - 1}A^k \mathbf{C} \left(\omega \circ  \phi^{-k}\right)
    \end{align*}
    is injective.
\end{lemma}

\noindent\textbf{Proof of the Lemma.\ \ }   Let $m _1, m _2 \in M _P $ be such that  $f_{(\phi, \omega,F)}(m_1) = f_{(\phi, \omega,F)}(m_2)$. This equality is equivalent to the following expressions:
\begin{align*}
\sum_{t = 0}^{\infty} A^t \mathbf{C} \omega \phi^{-t}(m_1) &= \sum_{t = 0}^{\infty} A^t \mathbf{C} \omega \phi^{-t}(m_2), \\
\sum_{t = 0}^{\infty} \sum_{k = 0}^{\ell-1} A^{t\ell + k} \mathbf{C} \omega \phi^{-(t\ell+k)}(m_1) &= \sum_{t = 0}^{\infty} \sum_{k = 0}^{\ell-1} A^{t\ell + k} \mathbf{C} \omega \phi^{-(t\ell+k)}(m_2), \\
\sum_{t = 0}^{\infty} (A^\ell)^t \sum_{k = 0}^{\ell-1} A^{k} \mathbf{C} \omega \phi^{-k}(m_1) &= \sum_{t = 0}^{\infty} (A^\ell)^t \sum_{k = 0}^{\ell-1} A^{k} \mathbf{C} \omega \phi^{-k}(m_2). 
\end{align*}
Given that $\rho(A)<1 $ then $\rho(A^{\ell})<1 $ necessarily and hence this equality can be rewritten as
\begin{equation*}
(I - A^\ell)^{-1} \sum_{k = 0}^{\ell-1} A^{k} \mathbf{C} \omega \phi^{-k}(m_1) =  (I - A^\ell)^{-1} \sum_{k = 0}^{\ell-1} A^{k} \mathbf{C} \omega \phi^{-k}(m_2),
\end{equation*}
which is equivalent to 
$\sum_{k = 0}^{\ell-1} A^{k} \mathbf{C} \omega \phi^{-k}(m_1) = \sum_{k = 0}^{\ell-1} A^{k} \mathbf{C} \omega \phi^{-k}(m_2) $ and hence, by definition, to $g_{(\phi, \omega,F)}(m_1) = g_{(\phi, \omega,F)}(m_2)$, which proves the statement. \quad $\blacktriangledown $

\medskip

If we now define  $\Phi_{(\ell, \omega)}: M \to \mathbb{R}^l$ as 
\begin{align*}
    \Phi_{(\ell, \omega)}(m) := \left( \omega(m) , \omega\circ  \phi^{-1}(m) , \ldots, \omega\circ\phi^{-(\ell-1)}(m) \right)^{\top},
\end{align*}
we note that the map $g_{(\phi, \omega,F)}  $  can be rewritten as $g_{(\phi, \omega,F)} = Q \Phi_{(\ell, \omega)}  $,  where $Q \in \mathbb{M}_{N, \ell} $ is a matrix whose $(k+1) $th-column is set to the vector $A ^k \mathbf{C} $. The hypotheses on the vectors $\left\{A ^j \mathbf{C}\right\}_{j \in \left\{0,1, \ldots, N-1 \right\}}$ forming a linearly independent set and that $N> \ell $ guarantee that ${\rm rank}\, Q = N $ and hence that the associated linear map $Q: \mathbb{R}^{\ell}\longrightarrow \mathbb{R}^N$ is injective. With this notation we now show that if $g_{(\phi, \omega,F)}  $ is not injective in $M _P $ then a perturbation $\omega ' \in C ^2(M, \mathbb{R}^N) $ of $\omega  $ can be chosen so that $g_{(\phi, \omega',F)}  $ is. More specifically, define
\begin{equation}
\label{perturbation for mp}
\omega':= \omega+\sum_{i=1}^P\sum_{j=1}^{n _i}\varepsilon_{ij}\Psi _{ij},
\end{equation}
where $\Psi _{ij}$ are bump functions with non-intersecting supports $U _{ij} $ such that $m _{ij}:=\phi^{-(j-1)}(m _i) \in U _{ij} $ and, moreover, $\Psi _{ij} \left(\phi^{-(j-1)}(m _i)\right)= \Psi _{ij}(m _{ij})=1/ \mathcal{L}(i,j)$. The symbol $\mathcal{L}(i,j) \in \mathbb{N} $  denotes the ordinal of the pair $(i,j) $ in lexicographic order.

We now show that the constants $\varepsilon_{ij} $ can be chosen so that $\omega'  $ is as close as we want to $\omega  $ and, at the same time,  $g_{(\phi, \omega',F)}  $ is injective. Firstly, it is easy to see that, by construction,
\begin{align*}
    \Phi_{(\ell, \omega')}(m) :=\Phi_{(\ell, \omega)}(m) +\sum_{i=1}^P\sum_{j=1}^{n _i}\varepsilon_{ij} \left( \Psi _{ij}(m) , \Psi _{ij}\circ  \phi^{-1}(m) , \ldots, \Psi _{ij}\circ\phi^{-(\ell-1)}(m) \right)^{\top}.
\end{align*}
Second, if $m_{i _1j _1} $ and $m_{i _2j _2} $ are two different periodic points then
\begin{equation}
\label{injectivity condition before lexi}
g_{(\phi, \omega',F)} (m_{i _1j _1})-g_{(\phi, \omega',F)} (m_{i _2j _2})=g_{(\phi, \omega,F)} (m_{i _1j _1})-g_{(\phi, \omega,F)} (m_{i _2j _2})+ Q\left(\varepsilon_{i _1 j _1}\mathbf{v}_{i _1 j _1}-\varepsilon_{i _2 j _2}\mathbf{v}_{i _2 j _2}\right) ,
\end{equation}
where the vectors $\mathbf{v}_{i _1 j _1} \in \mathbb{R}^{{\rm Card} M _P} $ have entries equal to zero except at the slots that are multiples of the period of the corresponding periodic point. More specifically, if the periodic point  $m_{i j }  $ has period $n _{ij}$, then
\begin{equation}
\label{definition of vs}
 (\mathbf{v}_{i  j }) _i:= \left\{
\begin{array}{l}
1/ \mathcal{L}(i,j) \quad \mbox{when $i=1$ or $i-1 $ is a multiple of $n_{ij}$,}\\
0 \quad \mbox{otherwise.}
\end{array}
\right.
\end{equation}
Using  the injectivity of $Q$ and Lemma \ref{injectivity_lemma}, we now show that we can choose the perturbation constants $\varepsilon_{i  j }$   so that the restriction of $g_{(\phi, \omega',F)} $ to $M _P $ is injective. Let
\begin{equation}
\label{choice of vareps}
\varepsilon_{i  j }:= \epsilon \|g_{(\phi, \omega,F)} (m_{i j })\|, \quad \mbox{for some constant $\epsilon>0$.} 
\end{equation}
We now show that if $\epsilon>0$ is chosen so that 
\begin{multline}
\label{choice of epsil}
\epsilon \max_{(i_1,j_1), (i_2,j_2)} \left\{\left\| Q\left(\left\|g_{(\phi, \omega,F)} (m_{i_1 j_1) }\right\|\mathbf{v}_{i _1 j _1}-\left\|g_{(\phi, \omega,F)} (m_{i_2 j_2} )\right\|\mathbf{v}_{i _2 j _2}\right)  \right\|\right\}\\
<\min_{(i_1,j_1), (i_2,j_2)} \left\{\left\| g_{(\phi, \omega,F)} (m_{i_1 j_1} )-g_{(\phi, \omega,F)} (m_{i_2 j_2} ) \right\|   \right\}
\end{multline}
then the injectivity of $g_{(\phi, \omega',F)} \mid _{M _P} $ is guaranteed. Indeed, consider first the case of two distinct periodic points $m_{i _1j _1} $  and  $m_{i _2j _2}$ for which $g_{(\phi, \omega,F)} $ fails to be injective, that is, $g_{(\phi, \omega,F)} (m_{i _1j _1})= g_{(\phi, \omega,F)} (m_{i _2j _2}) $. In that case, by \eqref{injectivity condition before lexi} and \eqref{choice of vareps} we have that
\begin{equation}
\label{second step}
g_{(\phi, \omega',F)} (m_{i _1j _1})-g_{(\phi, \omega',F)} (m_{i _2j _2})=
\epsilon \left\|g_{(\phi, \omega,F)} (m_{i _1j _1})\right\|Q\left(\mathbf{v}_{i _1 j _1}-\mathbf{v}_{i _2 j _2}\right).
\end{equation}
Given that $\mathbf{v}_{i _1 j _1}-\mathbf{v}_{i _2 j _2} \neq {\bf 0}$ (notice, for instance, that $\left(\mathbf{v}_{i _1 j _1}-\mathbf{v}_{i _2 j _2}\right)_1=1/ \mathcal{L}(i _1, j _1)- 1/ \mathcal{L}(i _2, j _2)\neq 0$) and $Q$ is injective then  $Q \left(\mathbf{v}_{i _1 j _1}-\mathbf{v}_{i _2 j _2}\right)\neq {\bf 0}$ and hence $g_{(\phi, \omega',F)} (m_{i _1j _1})\neq g_{(\phi, \omega',F)} (m_{i _2j _2}) $ necessarily by \eqref{second step}. In the case $g_{(\phi, \omega,F)} (m_{i _1j _1})\neq g_{(\phi, \omega,F)} (m_{i _2j _2}) $ the same conclusion can be drawn because the choice of $\epsilon>0$ in \eqref{choice of epsil} guarantees  that
\begin{equation*}
\left\|Q\left(\varepsilon_{i _1 j _1}\mathbf{v}_{i _1 j _1}-\varepsilon_{i _2 j _2}\mathbf{v}_{i _2 j _2}\right)\right\|<   \left\|g_{(\phi, \omega,F)} (m_{i _1j _1})-g_{(\phi, \omega,F)} (m_{i _2j _2})\right\|
\end{equation*}
which by \eqref{injectivity condition before lexi} ensures that, again, $g_{(\phi, \omega',F)} (m_{i _1j _1})\neq g_{(\phi, \omega',F)} (m_{i _2j _2}) $, as required.

We now show that if $\left. f_{(\phi, \omega,F)}\right|_{M _P}$ is injective then there exists an open set $V_P$ such that $M _P \subset V _P  $ and $\left. f_{(\phi, \omega,F)}\right|_{\overline{V _P}}$ is also injective. By the Immersion Theorem \cite[Theorem 3.5.7]{mta} we know that there exists $n \in \mathbb{N}  $ such that the balls $B_{2^{-n}}(m _{ij})$ do not intersect and that the restriction of $f_{(\phi, \omega,F)} $ to each of them is a collection of injective maps. It could still be, however, that the images of different balls intersect. The continuity of $f_{(\phi, \omega,F)} $ and the fact that $\left. f_{(\phi, \omega,F)}\right|_{M _P}$ is injective implies that $n$ can be chosen sufficiently high so that this does not happen. Indeed, if this was not the case for the balls around the periodic points, say, $m_{i _1j _1} $ and $m_{i _2j _2} $, then it would be possible to construct two sequences $\{m_{i _1j _1,l}\}_{l \in \mathbb{N}} $ and $\{m_{i _2j _2, l}\}_{l \in \mathbb{N}} $ with limits $m_{i _1j _1} $ and $m_{i _2j _2} $ for which $f_{(\phi, \omega,F)}(m_{i _1j _1,l})=f_{(\phi, \omega,F)}(m_{i _2j _2,l}) $ for each $ l \in \mathbb{N}  $. By continuity this implies that $f_{(\phi, \omega,F)}(m_{i _1j _1})=f_{(\phi, \omega,F)}(m_{i _2j _2}) $ which is in contradiction with the injectivity of $\left. f_{(\phi, \omega,F)}\right|_{M _P}$ and hence proves the injectivity of $ f_{(\phi, \omega,F)} $ restricted to $V _P = \bigcup_{ij}B_{2^{-n}}(m _{ij}) $, with $n$ chosen so that the properties of the corresponding balls designed above are satisfied. Notice that by doubling $n$, if necessary, it is also easy to ensure the injectivity of $\left. f_{(\phi, \omega,F)}\right|_{\overline{V _P}}$.

\medskip

\noindent {\bf Step 2: Global injectivity.} We firstly prove an important local intermediate result.

\begin{lemma}
\label{lemma1 with r}
If $M$ is a compact differentiable manifold endowed with a metric $d$ and $f : M \longrightarrow \mathbb{R}^N$ is an immersion, then there exists a constant $r > 0$ such that for any $m \in M$ the restriction $f\rvert_{B_r(m)}$ of $f$ to the open ball $B_r(m) \subset M$ of radius $r$ and center $m$ is injective.
\end{lemma}

\noindent\textbf{Proof.\ \ }The Immersion Theorem (\cite[Theorem 3.5.7]{mta}) implies that  each $m \in M$ has an open neighborhood $U_m \subset M$ such that $f\rvert_{U_m}$ is injective. The collection of sets $\{ U_m\}_{ m \in M }$ forms an open cover of $M$. Then, by Lebesgue's number lemma \cite[Lemma 27.5]{Munkres:topology}, there exists a $\delta > 0$ such that every set of diameter $\delta$ is contained in some set in the family $\{ U_m\}_{ m \in M }$. The lemma is proved by choosing  $r = \delta/2$. $\blacktriangledown $

\medskip

Since $M$ is compact and $f_{(\phi, \omega, F)} : M \longrightarrow \mathbb{R}^N$ is an immersion, this lemma implies the existence of a constant $r > 0$ such that for any $m \in M$ the restriction $f_{(\phi, \omega, F)}\rvert_{B_r(m)}$ of $f_{(\phi, \omega, F)}$ to the open ball $B_r(m)$ is  injective.
We now define the set $W \subset M \times M$ as follows using the open set $V_P$ whose existence we proved in Step 1.

\begin{equation*}
W := \{ (m_1 , m_2) \in (M \times M) \setminus (V_P \times V_P) \mid \ d(m_1,m_2) \geq r \}.
\end{equation*}
The set $W$ comprises pairs $(m_1,m_2) \in M$ whose entries satisfy one of two conditions:
\begin{enumerate}
\item Neither $m_1$ nor $m_2$ are in $V_P$.
\item One of $m_1$ and $m_2$ is in $V_P$ and the other is not.
\end{enumerate}
In view of this, the injectivity of $f_{(\phi, \omega, F)}\rvert_{V_P}$ proved in the Step 1 together with Lemma \ref{lemma1 with r} imply that if we show that $f_{(\phi, \omega, F)}(m _1)\neq f_{(\phi, \omega, F)}(m _2)$ for all $(m _1, m _2) \in W $ then $f_{(\phi, \omega, F)} $ is globally injective and the proof is concluded.

We start the proof of this fact by first defining, for each $m \in M$,  a collection of nested balls $\{ B_{2^{-n}}(m) \mid n \in \mathbb{N} \}$ centered at $m$ with radius $2^{-n}$. Let $(m_1,m_2) \in W$, and assume from now on without loss of generality that $m_1 \in W \setminus V_p$. Let $T(n,m_1,m_2)$ denote the largest integer such that the following two properties hold. Firstly, the sets
\begin{align*}
    \{B_{2^{-n}}(\phi^{-t}(m_1))\}_{t=0, \ldots, T(n,m_1,m_2)-1}
\end{align*}
are disjoint and secondly
\begin{align*}
    B_{2^{-n}}(\phi^{-t}(m_1)) \cap B_{2^{-n}}(\phi^{-s}(m_2)) = \emptyset \quad \mbox{for all} \quad t,s \in \{ 0 , \ldots , T(n,m_1,m_2) - 1\}.
\end{align*}
Notice now that by the continuity of $\phi$, for each $n \in \mathbb{N}$ and pair $(m_1,m_2) \in W$ there is an open neighbourhood $U_{(m_1,m_2)} \subset M \times M$ of $(m_1,m_2)$ such that $T(n,m_1',m_2') = T(n,m_1'',m_2'')$ for all $(m_1',m_2'),(m_1'',m_2'') \in U_{(m_1,m_2)}$. The collection $\{ U_{(m_1,m_2)} \ | \ (m_1,m_2) \in W \}$ covers $W$ and since it is a compact set we can extract a finite subcover $\{ U_a \ | \ a \in \mathcal{A} \}$, where $\mathcal{A}$ is a finite set. Then we can choose one pair $(m_1^a,m_2^a) \in U_a$ for each $a \in \mathcal{A}$ and notice that 
\begin{align*}
    \min_{(m_1,m_2) \in W} \{T(n,m_1,m_2)\} = \min_{\{(m_1^a,m_2^a) \ | \ a \in A\}} \{T(n,m_1^a,m_2^a)\}.
\end{align*}
The importance of this equality is that, since $\mathcal{A}$ is a finite set, the minimum on the right hand side is realized by a pair $(m_1^*,m_2^*) \in W$. Let $T(n) = T(n,m_1^*,m_2^*) = \min_{(m_1,m_2) \in W} T(n,m_1,m_2)$.
Observe that as $n \to \infty$ the families $\{ B_{2^{-n}}(\phi^{-t}(m_1^*)) \}_{t \in \mathbb{N}}$ and $\{ B_{2^{-n}}(\phi^{-t}(m_2^*)) \}_{t \in \mathbb{N}}$ converge to $\{\phi^{-t}(m_1^*)\}_{t \in \mathbb{N}}$ and $\{\phi^{-t}(m_2^*)\}_{t \in \mathbb{N}}$ respectively. The point $m_1^*$ is not periodic so the infinite orbit $\{\phi^{-t}(m_1^*)\}_{t \in \mathbb{N}}$ of singletons is disjoint, and furthermore does not intersect any point in $\{\phi^{-t}(m_2^*)\}_{t \in \mathbb{N}}$. This allows us to conclude that $T(n) \to \infty$ as $n \to \infty$.

The fact that we just proved guarantees the existence of a $\nu \in \mathbb{N}$ such that $T(\nu) = N$. Thus for all pairs $(m_1,m_2) \in W$, the collection
\begin{align*}
    \{B_{2^{-\nu}}(\phi^{-t}(m_1))\}_{t=0, \ldots, N-1}
\end{align*}
is disjoint and 
\begin{align*}
    B_{2^{-n}}(\phi^{-t}(m_1)) \cap B_{2^{-n}}(\phi^{-s}(m_2)) = \emptyset \quad \mbox{for all} \quad t,s \in \{ 0 , \ldots , N - 1\}.
\end{align*}
Now for any $n > \nu$ the collection
\begin{align*}
    \mathcal{C}_n = \{ B_{2^{-(n+1)}}(m) \ | \ m \in M \}
\end{align*}
forms on open cover of $M$ from which we can extract a finite subcover $\{B_i \ | \ i \in J_n\}$ for $J$ a finite set with cardinality  $\ell(n) \in \mathbb{N}$. Now define a partition of unity $\{ \lambda_i \ | \ i \in J_n \}$ subordinate to $\{B_i \ | \ i \in J_n\}$. We impose on this partition of unity the special property that for each $m \in M$ there exists an $i \in J_n$ such that $\lambda_i(m) \geq 1/2$. Now we define the perturbed observation function
\begin{align*}
    \omega_n = \omega + \sum^{\ell(n)}_{i=1} \epsilon_i \lambda_i
\end{align*}
where $\epsilon_i \in \mathbb{R}$ is the $i$th component of a vector $\boldsymbol{\epsilon} \in \mathbb{R}^{\ell(n)}$ with positive entries. Then we define $\Psi_n : M \times M \times \mathbb{R}^{\ell(n)} \to \mathbb{R}^N$ by
\begin{equation}
\label{perturbation lambda}
\Psi_n(m_1,m_2,\epsilon) = f_{(\phi, \omega_n, F)}(m_1) - f_{(\phi, \omega_n, F)}(m_2).
\end{equation}
Let $\Delta = \{ (m,m) \in M\times M \mid m \in M \}$ be the diagonal set. Given an arbitrary open neighborhood $\mathcal{N} \subset C^1(M,\mathbb{R})$ of the observation function $\omega \in C^2(M,\mathbb{R})$
our goal is to find $\boldsymbol{\epsilon} \in \mathbb{R}^{\ell(n)}$ such that $\omega_n \in \mathcal{N}$ and that for all $(m_1,m_2) \in (M \times M) \setminus \Delta$ we have that $\Psi_n(m_1,m_2,\boldsymbol{\epsilon}) \neq {\bf 0}$.

First of all, we observe that for any pair $(m_1,m_2) \in (M \times M) \setminus W$ either $d(m_1,m_2) < r$ or both $m_1,m_2 \in V_P$. In the former case, $\Psi_n(m_1,m_2,{\bf 0}) \neq  {\bf 0}$ unless $(m_1,m_2) \in \Delta$ by Lemma \ref{lemma1 with r}, and in the latter case, $\Psi_n(m_1,m_2,{\bf 0}) \neq {\bf 0}$ unless $(m_1,m_2) \in \Delta$ because $f_{(\phi, \omega, F)}\rvert_{V_P}$ is injective by the Step 1. Now $\Psi_n$ is continuous so there is an open neighbourhood $U_{{\bf 0}} \subset \mathbb{R}^{\ell(n)}$ of ${\bf 0} \in \mathbb{R}^{\ell(n)}$ such that for all $\boldsymbol{\epsilon} \in U_{{\bf 0}}$ we have $\Psi_n(m_1,m_2,\boldsymbol{\epsilon}) \neq {\bf 0}$ for all $(m_1,m_2) \in (M \times M) \setminus W$ unless $(m_1,m_2) \in \Delta$. So all that remains is to find  $\boldsymbol{\epsilon} \in U_{{\bf 0}} \subset \mathbb{R}^{\ell(n)}$ such that $\Psi_n(m_1,m_2,\boldsymbol{\epsilon}) \neq 0$ for all $(m_1,m_2) \in W$.

We start by noting that if ${\bf 0} \in \mathbb{R} ^N $  is not in the image of $\Psi_n \rvert_{W \times \{{\bf 0}\}}$ then we are done so we shall assume the opposite. In that case we proceed by showing that  $\Psi_n \rvert_{W \times \{{\bf 0}\}}$ is a submersion. If that is the case, then for some open set $X \subset (M \times M \times \mathbb{R}^{\ell(n)})$ containing $W \times \{{\bf 0}\}$ then the restriction $\Psi_n \rvert_{X}$ is also a submersion and hence by the Submersion Theorem \cite[Theorem3.5.4]{mta} the inverse image $\Psi_n\lvert_{X}^{-1}\left({\bf 0}\right)$ is a closed submanifold of dimension $2q+\ell(n) - N$ of the open submanifold $X\subset M \times  M \times \mathbb{R}^{\ell}$.
Moreover, if  $\pi : M \times M \times \mathbb{R}^{\ell(n)} \to \mathbb{R}^{\ell(n)}$ is the canonical projection defined by $\pi(m_1,m_2,\boldsymbol{\epsilon}) := \boldsymbol{\epsilon}$,  in these circumstances the complement $\mathbb{R}^{\ell(n)} \setminus \pi\left(\Psi_n\rvert_{X}^{-1}({\bf 0})\right)$ is a dense subset of $\mathbb{R}^{\ell(n)}$. Indeed, since $\pi$ is a continuously differentiable map, then so is its restriction  $\pi\rvert_{\Psi_n\lvert_{X}^{-1}\left({\bf 0}\right)\times \mathbb{R}^{\ell}}: \Psi_n\lvert_{X}^{-1}\left({\bf 0}\right)\longrightarrow \mathbb{R}^{\ell} $ which by \cite[Chapter 3, Proposition 1.2]{Hirsch:book} guarantees the density of $\mathbb{R}^{\ell(n)} \setminus \pi\left(\Psi_n\rvert_{X}^{-1}({\bf 0})\right)$.  This implies that we can choose $\boldsymbol{\epsilon} \in \left(\mathbb{R}^{\ell(n)} \setminus \pi \left(\Psi_n\rvert_{X}^{-1}({\bf 0})\right)\right)$ as small as we want so that $\boldsymbol{\epsilon} \in U_{{\bf 0}}$ and $\omega_n \in \mathcal{N}$. We fix this $\epsilon$ and see that for any $(m_1,m_2) \in W$ the map $\Psi_n(m_1,m_2,\boldsymbol{\epsilon}) \neq  0 $, as required. 
Consequently, all that remains to be done is to find $n$ sufficiently large so that $\Psi_n \rvert_{W \times \{{\bf 0}\}}$ is a submersion, and then the proof will be complete.

We start by observing that by \eqref{perturbation lambda} 
\begin{equation*}
\omega_n \circ \phi^{-t} = \omega \circ \phi^{-t} + \sum_{i = 1}^{\ell(n)} \epsilon_i \lambda_i \circ \phi^{-t}
\end{equation*}
hence 
\begin{equation*}
\frac{\partial (\omega_n\circ  \phi^{-t})}{\partial \epsilon_j} = \lambda_j \circ \phi^{-t}.
\end{equation*}
Now we consider an arbitrary $(m_1,m_2) \in W$ assuming once again without loss of generality that $m_1 \in M \setminus V_p.$ For each point in the orbit $\{ \phi^{-t}(m_1) \}_{t = 0, \ldots , T(n)-1}$ there exists a $j(t) \in J_n$ such that $\lambda_{j(t)}\left( \phi^{-t}(m_1) \right)\geq 1/2$ by the special property that we imposed earlier  on the partition of unity $\{ \lambda_i \ | \ i \in J _n \}$. Now the support of $\lambda_{j(t)}$ is a ball $B_{j(t)}$ of radius $2^{-(n+1)}$ which contains $\phi^{-t}(m_1)$. Hence the ball $B_{j(t)} \subset B_{2^{-n}}(\phi^{-t}(m_1))$. Now since the sets in the family $\{B_{2^{-n}}(\phi^{-t}(m_1))\}_{t = 0 , \ldots , T(n) - 1}$ are disjoint then so are $\{ B_{j(t)} \}_{t = 0 , \ldots , T(n) - 1}$. 
Furthermore, since $B_{2^{-n}}(\phi^{-t}(m_1)) \cap B_{2^{-n}}(\phi^{-s}(m_2)) = \emptyset$ for all $t,s \in \{ 0 , \ldots , T(n) - 1\}$ hence $\lambda_{j(t)}(\phi^{-t}(m_2)) = 0$ for $t \in \{ 0, \ldots, T(n) - 1\}$. Thus
\begin{align*}
\frac{\partial (\omega_n\circ  \phi^{-t})}{\partial \epsilon_{j(t)}}(m_2) = 0.
\end{align*}
Now,
\begin{align*}
\Psi_n(m_1,m_2,\epsilon) &= \sum_{t = 0}^{T(n)-1} A^t \mathbf{C} (\omega_n (\phi^{-t}(m_1)) - \omega_n (\phi^{-t}(m_2)) ) \\ &+
\sum_{t = T(n)}^{\infty} A^t \mathbf{C} (\omega_n( \phi^{-t}(m_1)) - \omega_n( \phi^{-t}(m_2) ))
\end{align*}
hence for $t = 0, \ldots, T(n) - 1$
\begin{align*}
\frac{\partial \Psi_n}{\partial \epsilon_{j(t)}}(m_1,m_2,\epsilon) &= A^t \mathbf{C} (\lambda_{j(t)} (\phi^{-t}(m_1))) \\ &+
\sum_{t = T(n)}^{\infty} A^t \mathbf{C} (\lambda_{j(t)} (\phi^{-t}(m_1)) - \lambda_{j(t)}( \phi^{-t}(m_2) ) ). 
\end{align*}
By assumption $\{A^t \mathbf{C}\}_{t = 0, ... , N-1}$ are linearly independent, hence the vectors $\{A^t \mathbf{C} (\lambda_{j(t)} (\phi^{-t}(m_1)))\}_{t \in \{ 0 , \ldots, T(n)-1\}}$ necessarily span $\mathbb{R}^N$ because since $n> \nu  $  then $T(n)\geq N $. Crucially, for any $n$ the property $\lambda_{j(t)} (\phi^{-t}(m_1)) \geq 1/2$ holds and therefore, the residual term
\begin{align*}
\sum_{t = T(n)}^{\infty} A^t \mathbf{C} (\lambda_{j(t)} \phi^{-t}(m_1) - \lambda_{j(t)} \phi^{-t}(m_2) )
\end{align*}
may only spoil the spanning property of the vectors $\{A^t \mathbf{C} (\lambda_{j(t)} \phi^{-t}(m_1))\}_{t = 0, \ldots , N-1}$ if it is sufficiently large. Since by hypothesis $\rho(A)<1 $, the residual term converges uniformly over $(m_1,m_2) \in W$ to $0$ as $n$ grows. We choose consequently $n$ large enough so that for all $(m_1,m_2) \in W$ the residual term is too small to spoil the spanning property of $\{A^t \mathbf{C} (\lambda_{j(t)} \phi^{-t}(m_1))\}_{t \in \{ 0, \ldots , N-1\}}$. With this choice of $n$ we have that $\Psi_n\rvert_{W \times \{ 0 \}}$ is a submersion and the proof is complete. \quad $\blacksquare$

\paragraph{Linear reservoir embeddings} We conclude the theoretical part of the paper by showing that the embeddings whose existence we proved in Theorem \ref{the immersions are embeddings} using generic observation maps $\omega \in C ^2(M, \mathbb{R}^N)$ may be almost surely obtained, as it is customary in reservoir computing, by randomly drawing the connectivity matrix $A$ and the vector $\mathbf{C}$ of the linear system $F(\mathbf{x},z):=A \mathbf{x}+ \mathbf{C}z $. This result hinges on an important fact in random matrix theory whose proof has been kindly communicated to us by Friedrich Philipp and that is contained in the following statement. We recall that a random variable $X : \Omega \longrightarrow T $ defined on a probability space $(\mathbb{P}, \mathcal{F}, \mathbb{P}) $ and with values on a Borel measurable space $T $ is {\it regular} or {\it non-singular} whenever $\mathbb{P} \left(X =a\right)=0 $ for all $a \in T $.

\begin{proposition}[Friedrich Philipp]
\label{random_matrix_lemma}
Let $N \in \mathbb{N}  $,  $A\in \mathbb{M}_{N,N}$, and $\mathbf{C}\in \mathbb{R}^N$  and assume that the entries of $A$ and $\mathbf{C}$ are drawn using independent regular real-valued distributions. Then the following statements hold:
\begin{description}
\item [(i)]  The vectors $\mathbf{C}, A\mathbf{C}, A^2 \mathbf{C}, \ldots , A^{N-1}\mathbf{C}$ are linearly independent almost surely
\item [(ii)]  Given $m$ distinct complex numbers $\lambda_1, \ldots, \lambda_m \in \mathbb{C}$, where $m \leq N$, the event that $1, \lambda_1 , ... , \lambda_m \notin \sigma(A)$ ($\sigma(A)$ is the spectrum of $A$) and that the vectors
\begin{align*}
    (\mathbb{I} - \lambda_j A)^{-1}(\mathbb{I} - A)^{-1}(\mathbb{I} - A^N)\mathbf{C}, \quad \mbox{$j = 1, \dots , m$}
\end{align*}
are linearly independent holds almost surely.
\end{description}
\end{proposition}

\noindent\textbf{Proof.\ \ } The vectors $\mathbf{C}, A\mathbf{C}, A^2 \mathbf{C}, \ldots , A^{N-1}\mathbf{C}$ are linearly independent if and only if
\begin{align*}
    \text{det}\left(\mathbf{C}|A\mathbf{C}|A^2\mathbf{C}| \cdots | A^{N-1}\mathbf{C}\right) = 0
\end{align*} 
which, in the notation of Lemma \ref{polynomial_lemma}, can be written as
\begin{align*}
    \text{det}\left(p_0(A)\mathbf{C},p_1(A)\mathbf{C},p_2(A)\mathbf{C}, \ldots , p_{N-1}(A)\mathbf{C}\right) = 0
\end{align*}
 using the linearly independent polynomials $p_j(A) := A^j$, $j \in \left\{0, \ldots, N-1\right\}$. Part {\bf (i)} of the statement hence follows directly from Lemma \ref{polynomial_lemma}. Now we turn our attention to part {\bf (ii)}. First of all, $\lambda_j$ is an eigenvalue of $A$ if and only if $\lambda_j$ is a root of the characteristic polynomial of $A$. This event has probability $0$ by Lemma \ref{first lemma polynomial} and hence $1, \lambda_1 , \ldots , \lambda_m \notin \sigma(A)$ almost surely. On this event, the inverses $(\mathbb{I} - \lambda_j A)^{-1}$ and $(\mathbb{I} - A)^{-1}$ exist. Furthermore, the product
 \begin{align*}
     \prod_{i=1}^{m}(\mathbb{I} - \lambda_i A)
 \end{align*}
 is an invertible matrix. Therefore, the vectors 
 \begin{align*}
      (\mathbb{I} - \lambda_j A)^{-1}(\mathbb{I} - A)^{-1}(\mathbb{I} - A^N)\mathbf{C}, \quad \mbox{with  $j = 1, ... , m$,} 
 \end{align*}
are linearly independent if and only if 
 \begin{align}
\label{vecs}
     &\prod_{i=1}^{m}(\mathbb{I} - \lambda_i A)(\mathbb{I} - \lambda_j A)^{-1}(\mathbb{I} - A)^{-1}(\mathbb{I} - A^N)\mathbf{C},  \quad \mbox{with  $j = 1, ... , m$,} 
 \end{align}
are linearly independent.
We can now rewrite the vectors in \eqref{vecs} as
 \begin{equation*}
    \prod_{i=1}^{m}(\mathbb{I} - \lambda_i A)(\mathbb{I} - \lambda_j A)^{-1}(\mathbb{I} - A)^{-1}(\mathbb{I} - A^N)\mathbf{C} 
    =\prod_{i \neq j}^{m}(\mathbb{I} - \lambda_i A)(\mathbb{I} - A)^{-1}(\mathbb{I} - A^N)\mathbf{C} 
     = \prod_{i \neq j}^{m}(\mathbb{I} - \lambda_i A)\sum_{k=0}^{N-1} A^k \mathbf{C},
 \end{equation*}
where we used the relation 
\begin{equation*}
(\mathbb{I} - A^N)=(\mathbb{I}-A)\sum_{k=0}^{N-1}A ^k \quad \mbox{and hence that} \quad (\mathbb{I}-A)^{-1}(\mathbb{I} - A^N)=\sum_{k=0}^{N-1}A ^k.
\end{equation*}
Now, if we are able to show that the family
 \begin{align*}
     p_j(x) = \prod_{i \neq j}^{m}(1- \lambda_i x)\sum_{k=0}^{N-1} x^k , \quad \mbox{with } \quad j \in \left\{ 1, \ldots , m\right\}
 \end{align*}
is linearly independent, then we can conclude by Lemma \ref{polynomial_lemma} that the vectors \eqref{vecs} are linearly independent almost surely, which would complete the proof. This is indeed the case because if $\mu_1, \ldots, \mu_n \in \mathbb{R} $ are such that 
\begin{equation*}
\sum_{j=1}^n \mu _jp _j(x)=0 \quad \mbox{then} \quad  \left(\sum_{k=0}^{N-1} x^k\right)\left(\sum_{j=1}^n \mu _j\prod_{i \neq j}^{m}(1- \lambda_i x)\right)=0.
\end{equation*}
Given that the polynomial $\sum_{k=0}^{N-1} x^k $ is non-zero, the previous equality is equivalent to $\sum_{j=1}^n \mu _j\prod_{i \neq j}^{m}(1- \lambda_i x)=0 $ which, evaluated at $x=1/ \lambda _k  $, implies that
\begin{equation*}
0=\sum_{j=1}^n \mu _j\prod_{i \neq j}^{m}\left(1- \lambda_i \frac{1}{\lambda_k}\right)=\mu _k\prod_{i \neq k}^{m}\left(1- \lambda_i \frac{1}{\lambda_k}\right).
\end{equation*}
Given that, by hypothesis, the values $\lambda_1, \ldots, \lambda_m \in \mathbb{C}$ are all different, we can conclude that $\prod_{i \neq k}^{m}\left(1-  \frac{\lambda_i}{\lambda_k}\right)\neq 0 $ and hence $\mu_k=0 $, necessarily. Since procedure can be repeated to obtain that $\mu_1, \ldots, \mu_n =0$, the result follows. \quad $\blacksquare$

\medskip

This proposition together with Theorem \ref{the immersions are embeddings} can be used to prove the following statement which is the main result of the paper.

\begin{theorem}[Linear reservoir embeddings]
\label{Theorem Linear reservoir embeddings}
Let $\phi \in {\rm Diff}^2(M)$ be a dynamical system on a compact manifold $M$ of dimension $q$ that exhibits finitely many periodic orbits. Suppose that for each periodic orbit $m$ of $\phi $ with period $n \in \mathbb{N} $, the derivative $T_m \phi^{-n} $ has $q$ distinct eigenvalues $\lambda_1, \lambda _2, \ldots, \lambda _q $. Let now $\ell \in \mathbb{N} $ be the lowest common multiple of all the periods of the finite periodic points of $\phi$ and let $N \in \mathbb{N} $ such that $N>\max \left\{2q, \ell\right\}$.

Construct now $\overline{A}\in \mathbb{M}_{N,N}$ and $\overline{\mathbf{C}}\in \mathbb{R}^N$  by drawing their entries using independent regular real-valued distributions. Then, there exist rescaled versions $A $ and $\mathbf{C}$ of $\overline{A}$ and $\overline{\mathbf{C}} $  respectively such that  the generalized synchronization $f_{(\phi, \omega,F)} \in C ^2(M, \mathbb{R}) $ associated to the state map $F(\mathbf{x},z ):=A \mathbf{x}+ \mathbf{C}z  $ is almost surely an embedding for generic $\omega \in C ^2(M, \mathbb{R}^N)$.
\end{theorem}

\noindent\textbf{Proof.\ \ } Proposition \ref{random_matrix_lemma} guarantees that the randomly drawn elements $\overline{A}$ and $\overline{\mathbf{C}} $ satisfy almost surely the hypotheses in parts {\bf (i)} and {\bf (ii)} of the statement of Theorem \ref{Theorem immersion}. However, in order to be able to invoke Theorem \ref{the immersions are embeddings}, we need to use a linear state map $F$ whose connectivity matrix $A$ is such that $\rho(A)<1 $ and, for any observation map $\omega \in C ^2(M, \mathbb{R})$, the corresponding generalized synchronization  $f_{(\phi, \omega,F)} \in C ^2(M, \mathbb{R}^N) $ and the map $\Theta_{(\phi, F)}: C ^2(M, \mathbb{R})  \longrightarrow  C ^2(M, \mathbb{R}^N)$ introduced in \eqref{functor gs} are continuous. It is obvious from parts {\bf (i)} and {\bf (ii)} in Proposition \ref{generalized synch with spectral radius} that this can be achieved by rescaling the matrix $\overline{A}$ and hence the statement follows from Theorem \ref{the immersions are embeddings}. \quad $\blacksquare$

\section{Numerical illustrations of attractor reconstruction, filtering, and forecasting}
\label{Numerical illustrations}

In this section we illustrate how the embeddings proposed in Theorem \ref{Theorem Linear reservoir embeddings} are able to reconstruct the attractor of various dynamical systems out of one-dimensional observations and, additionally, we show that randomly generated linear reservoir systems are  {efficient in the filtering and prediction of  dynamical systems observations in the presence of additive noise}.

Following the prescription proposed in the statement of Theorem \ref{Theorem Linear reservoir embeddings}, we shall randomly generate linear systems of the form $F(\mathbf{x},z)=A \mathbf{x}+ \mathbf{C}z $  to which we shall feed in the input variable $z$ finite-length one-dimensional $\omega$ observations  of three different dynamical systems $\phi$, namely, 
the Rossler system,  the Van der Pol oscillator, and the Lorenz system. For each of these systems we shall create reservoir states $\mathbf{x} _t $ according to the recursion
\begin{align}
\label{state generation algorithm}
    \mathbf{x}_{t+1} = A\mathbf{x}_t + C\omega (\phi^{t}(m)).
\end{align}
Due to the results  in the paper, we expect that the states $\mathbf{x}_t$ shall approximate  $f_{(\phi,\omega,F)} ( \phi^t(m))$ as $t \to \infty$ where $f_{(\phi,\omega,F)}$ is the corresponding embedding generalized synchronization introduced in Section \ref{Definitions and preliminary discussion}. This will be done in practice by keeping only the states $\mathbf{x} _t $  for all $t > T$, for some $T\in  \mathbb{N}$ where $[0,T]$ is called the washout period. The embedding properties of $f_{(\phi,\omega,F)}$ will become apparent in plots that will show that the dynamics of the original system and the state dynamics induced by its observations are topologically conjugate.

In order to illustrate the embedding properties of $f_{(\phi,\omega,F)}$ for the three dynamical systems we set up a reservoir system by following the steps:
\begin{enumerate}
    \item Randomly generate a 7 by 7 matrix $A'$ with IID uniform entries in the interval $[-0.5,0.5]$.
    \item Define the reservoir matrix $A := A' / \lVert A' \rVert$.
    \item Randomly generate a vector $\mathbf{C}\in \mathbb{R} ^7$  with IID uniform entries in the interval $[-0.5,0.5]$.
\end{enumerate}

\begin{figure}
  \centering
    \includegraphics[width=0.8\textwidth]{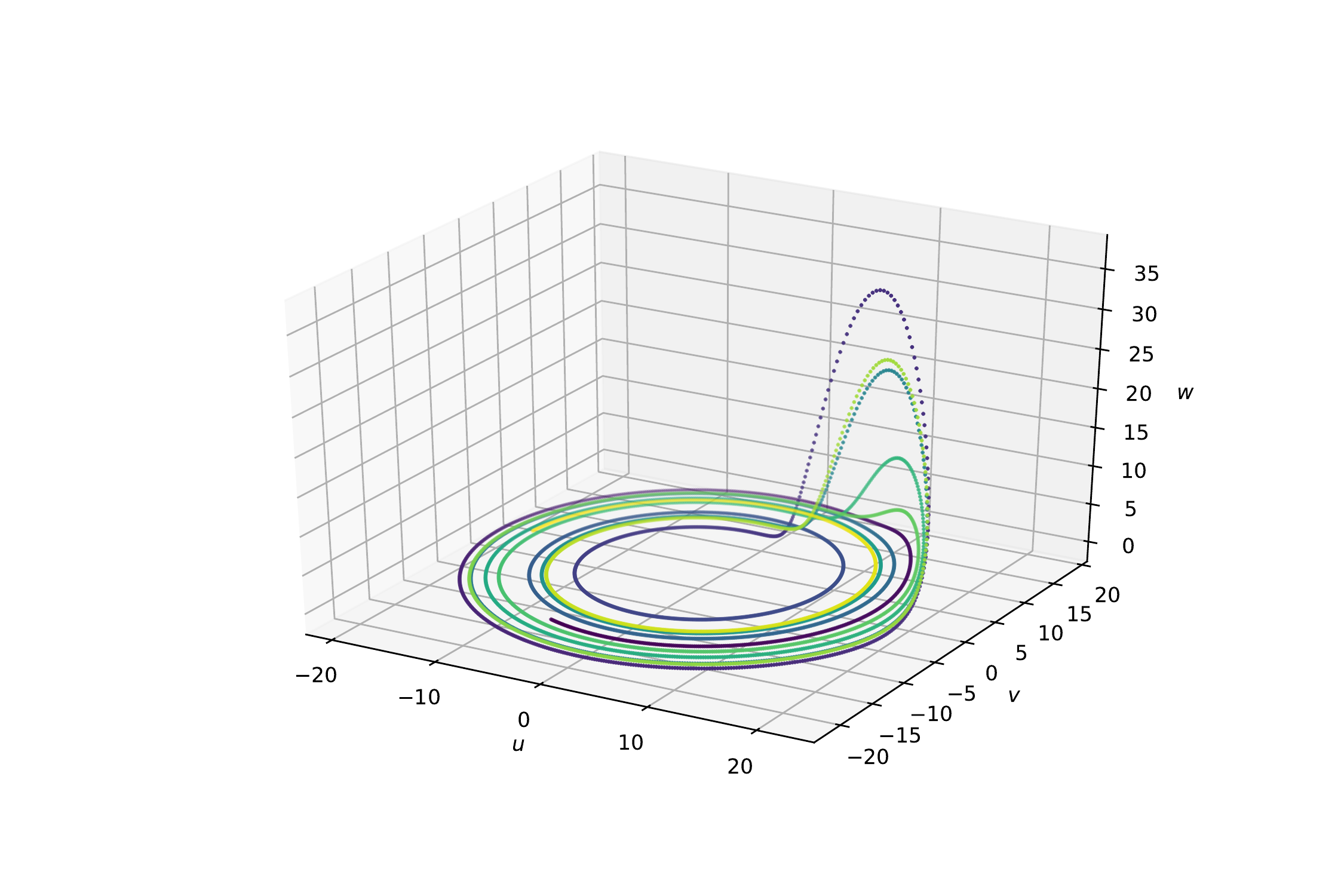}
 \caption{A trajectory of the R\"ossler system plotted for times in the interval $(60,120)$. Points at the start of the trajectory are purple, and points later on are yellow.}
 \label{fig::PRS_RS}
\end{figure}

\paragraph{The R\"ossler System.}
The R\"ossler system under a popular choice of parameters is described by the differential equations:
\begin{align*}
    \dot{u} &= -v -w,\\
    \dot{v} &=  u + v/10,\\
    \dot{w} &= 1/10 + w(u - 14).
\end{align*}
Using Python 3.7 and \texttt{scipy.integrate.odeint} we simulate a trajectory of the R\"ossler system
from the initial condition $(u_0,v_0,w_0) = (2,1,5) $
for $T=120$ time units, with time step $h = 0.01$. The result is plotted in Figure \ref{fig::PRS_RS} after using the interval $[0,60]$ as washout period. 



\begin{figure}[h!]
  \centering
    \includegraphics[width=0.8\textwidth]{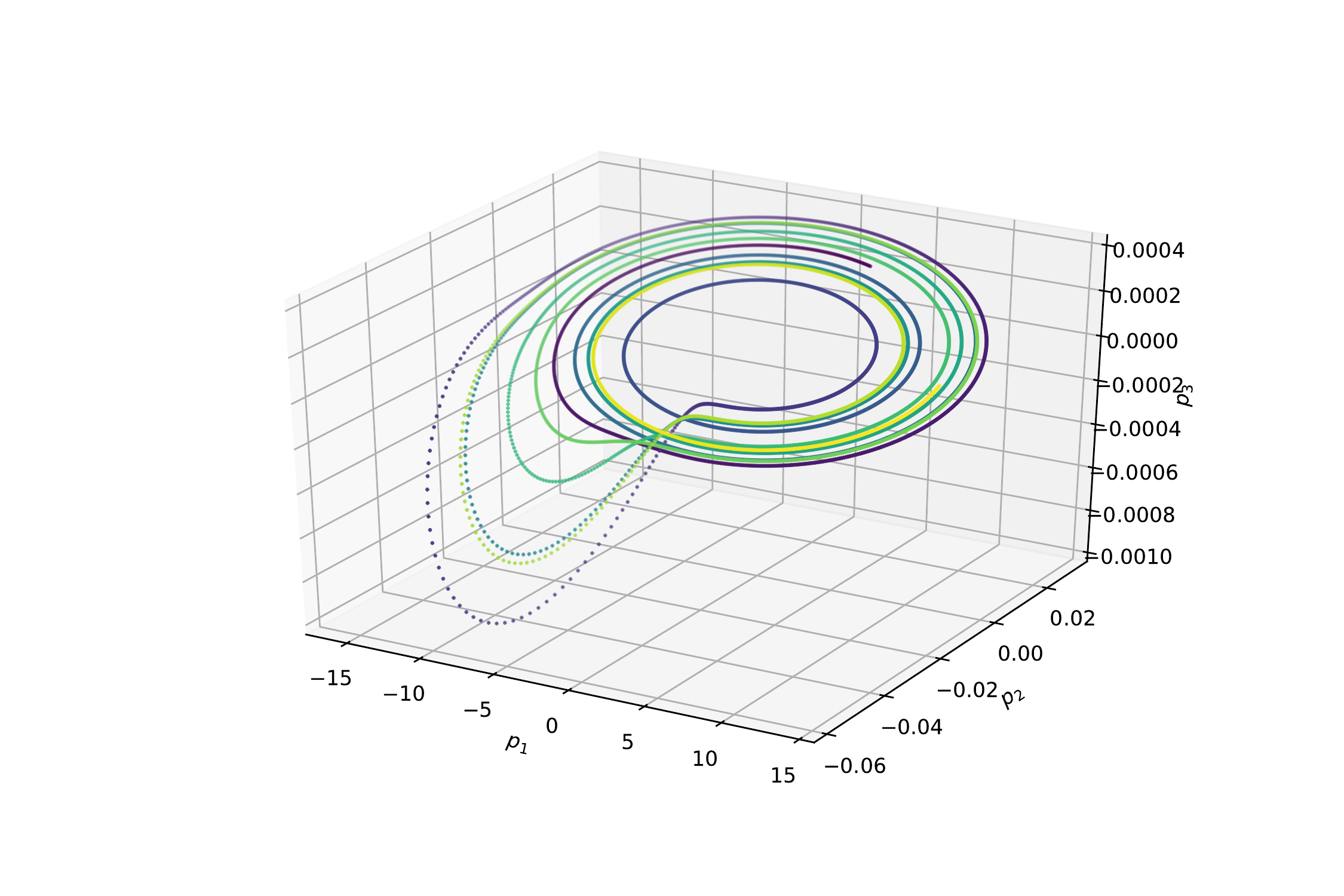}
\caption{Projection onto the first three principal components of the reservoir states in the interval $(60,120)$ for a system driven by the $u$-components of the R\"ossler system. Points at the start of the trajectory are purple, and points later on are yellow. It is worth pointing out the impressive resemblance with the original dynamical system even though this picture has been generated only using one of its components.}
    \label{fig::PRS_RS_reservoir}
\end{figure}

We then generate reservoir states using the recursion \eqref{state generation algorithm} and taking as observation map $\omega(u,v,w):=u $, that is, the first component of the R\"ossler system. 
A depiction of the projection of the corresponding states $x_t$ onto the first 3 principal components  is shown in Figure \ref{fig::PRS_RS_reservoir}.

%


\begin{figure}[h!]
  \centering
    \includegraphics[width=0.7\textwidth]{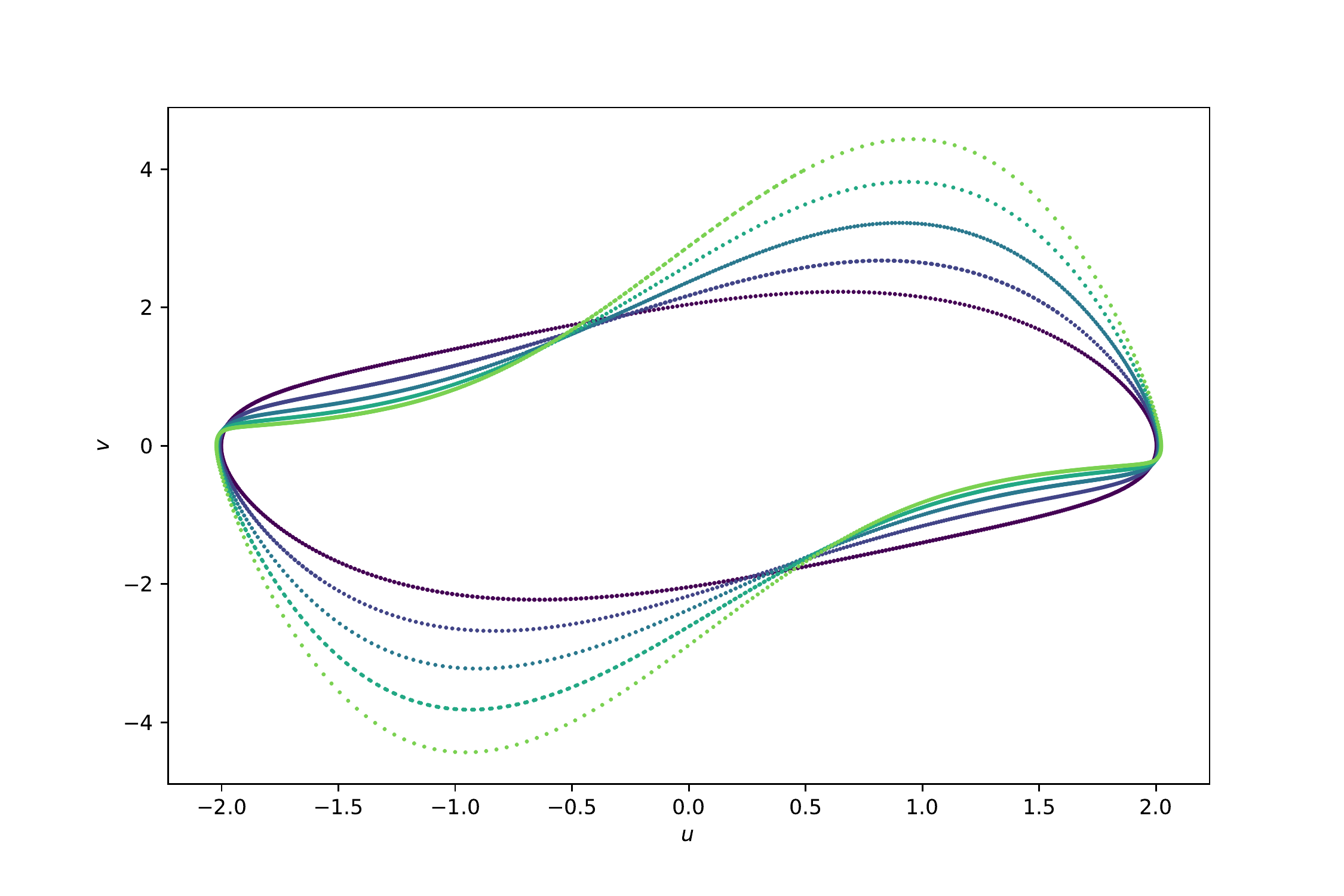}
  \caption{Limit cycles of the Van der Pol oscillator using damping parameters $\mu = 0.5,1,1.5,2,2.5$ plotted for times in the interval $(30,40)$. Darker colors denote lower values of $\mu$, brighter colors denote higher values of $\mu$.}
   \label{fig::PRS_VDP}
\end{figure}

\begin{figure}[h!]
  \centering
    \includegraphics[width=0.7\textwidth]{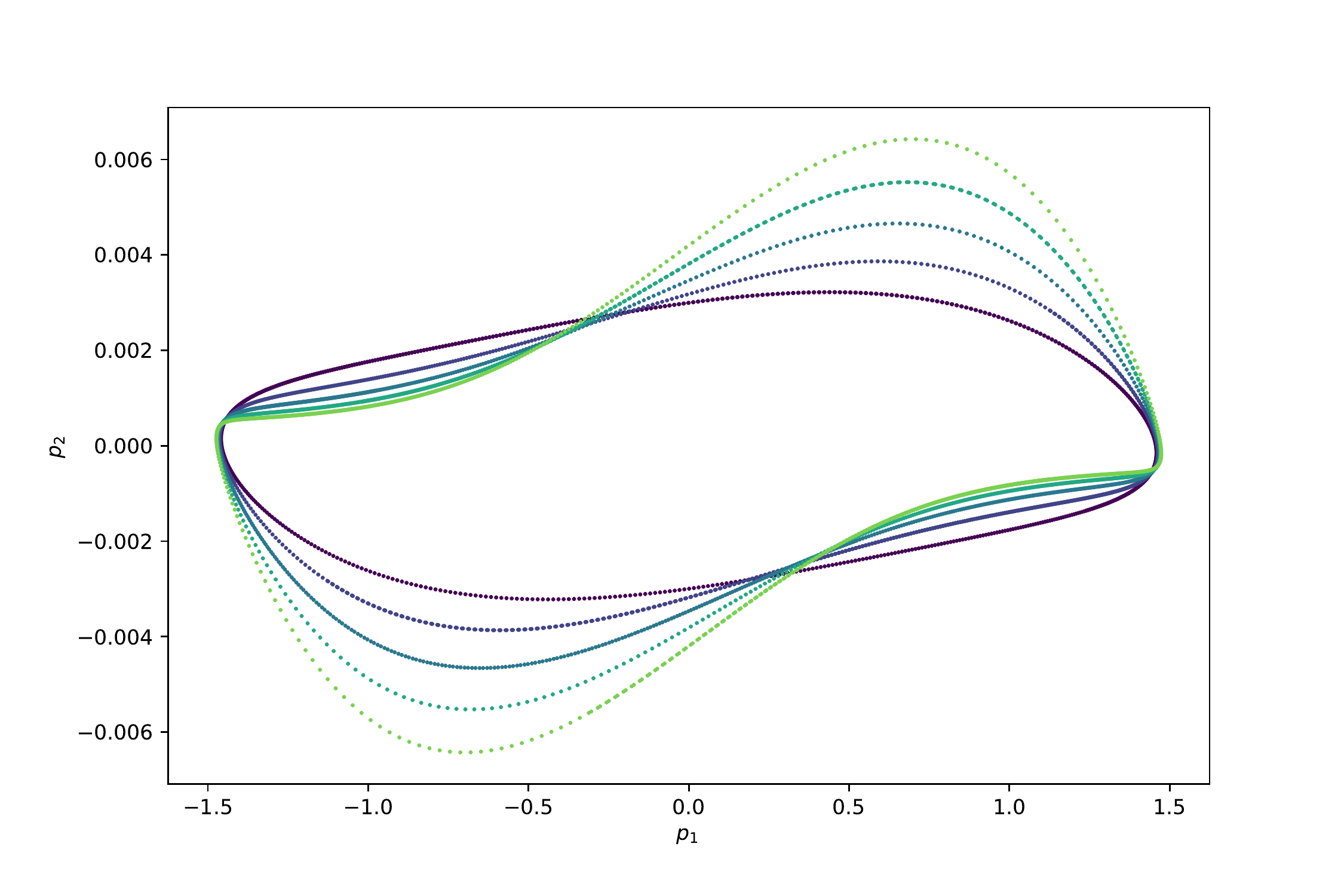}
  \caption{Projection onto the first two principal components of the reservoir states in the interval $(30,40)$ for a system driven by the $u$-components of the Van der Pol oscillator with damping parameters $\mu = 0.5,1,1.5,2,2.5$. Darker colors denote lower values of $\mu$, brighter colors denote higher values of $\mu$.}
  \label{fig::PRS_VDP_reservoir}
\end{figure}

\paragraph{The Van der Pol Oscillator.} We now repeat the same embedding procedure for the limit cycle of the Van der Pol oscillator with damping parameter $\mu$ is described by the differential equation in two dimensions
\begin{align*}
    \dot{u} &= v,\\
    \dot{v} &= \mu(1-u^2)v - u.
\end{align*}
Using the same discretization scheme as before we simulated trajectories of the Van der Pol oscillator for 40 time units with time step $h = 0.01$ for five different damping parameter values $\mu = 0.5, 1, 1.5, 2, 2.5$ but using always the same initial condition $(u_0,v_0) = (-4,5)$.
The result is plotted in Figure \ref{fig::PRS_VDP}.
Regarding the reservoir embedding we use a random linear system of dimension five using the same distributions for the entries as in previous paragraphs and we use as input the first component of the Van der Pol oscillator. A depiction of the projection of the corresponding states $x_t$ onto the first two principal components  is shown in Figure \ref{fig::PRS_VDP_reservoir} for each of the five different damping parameter values under consideration.

\begin{figure}[h!]
\centering
\includegraphics[width=.75\textwidth]{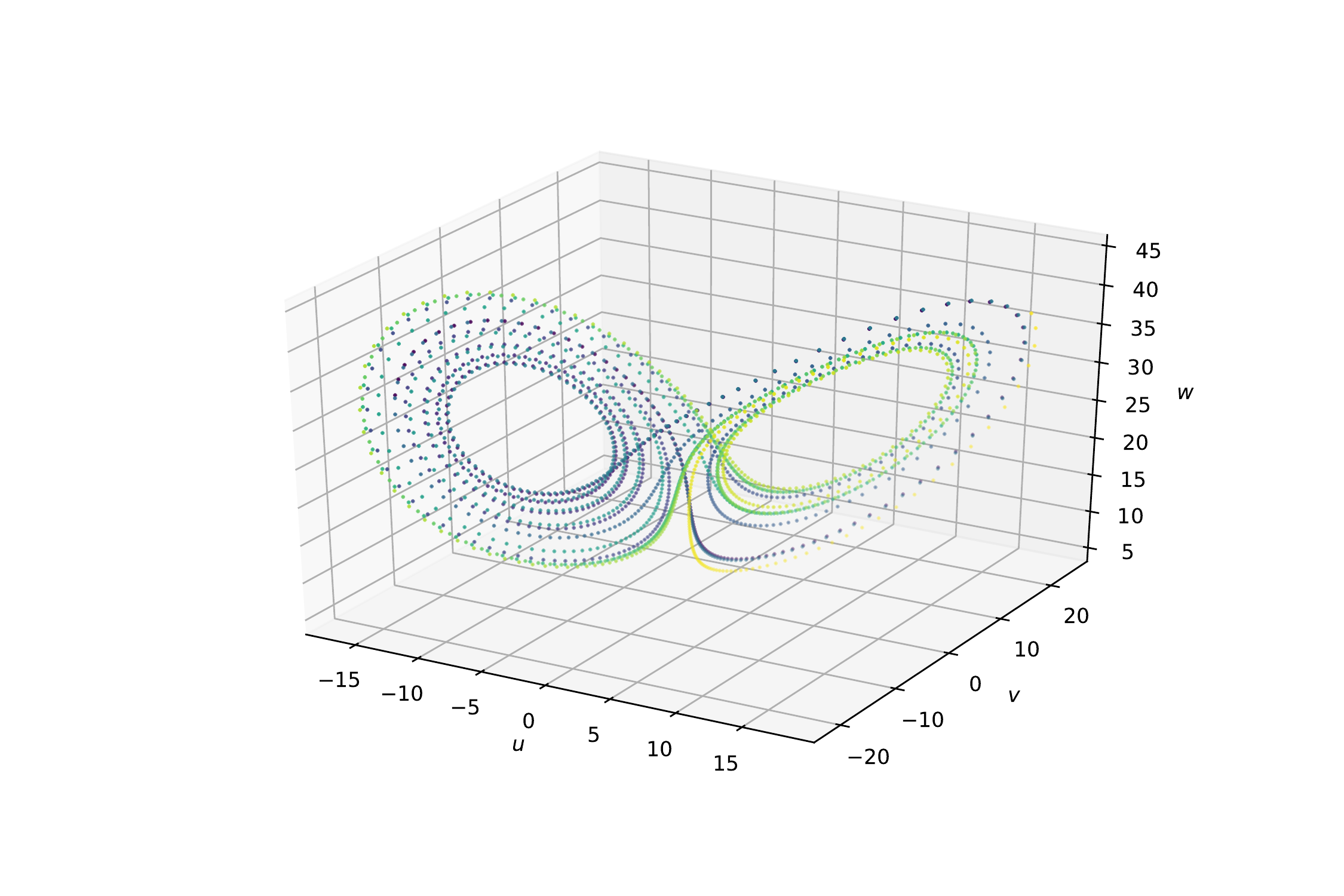}
\caption{A trajectory of the Lorenz system plotted for times in the interval $(20,40)$. Points at the start of the trajectory are purple, and points later on are yellow.}
\label{fig::PRS_lorenz}
\end{figure}

\begin{figure}[h!]
 \centering
    \includegraphics[width=0.8\textwidth]{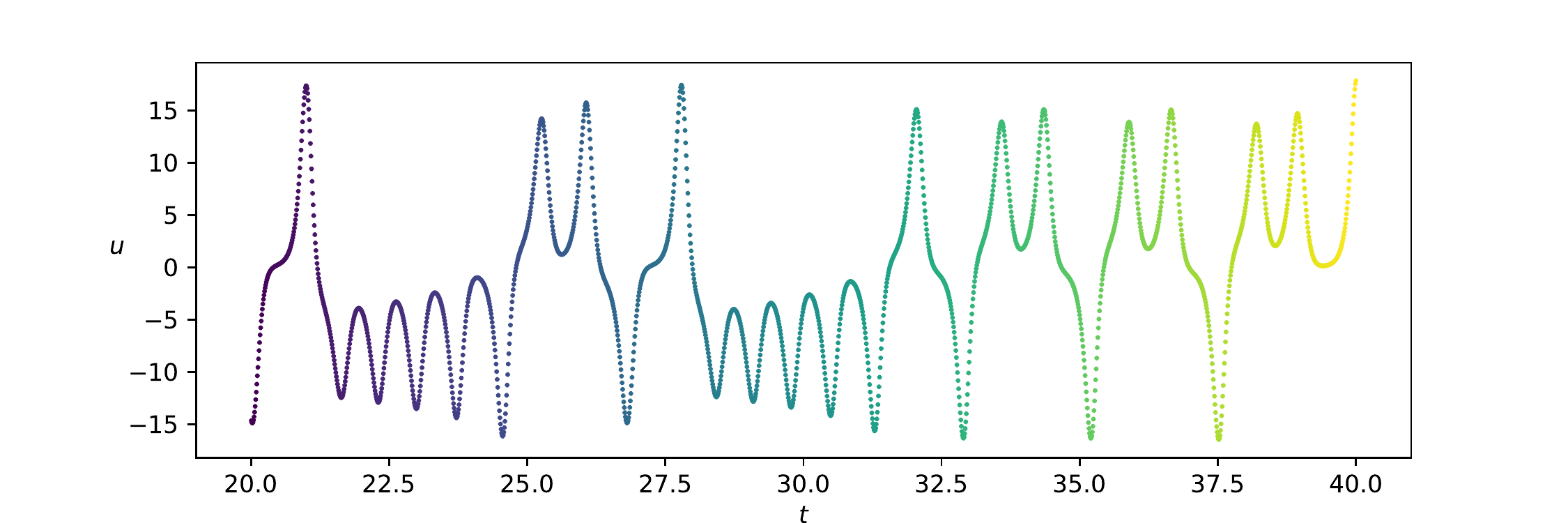}
\caption{A trajectory of the $u$-component of the Lorenz system plotted for times in the interval $(20,40)$. Points at the start of the trajectory are purple, and points later on are yellow.}
\label{fig::PRS_lorenz_x}
\end{figure}

\begin{figure}[h!]
  \centering
    \includegraphics[width=0.8\textwidth]{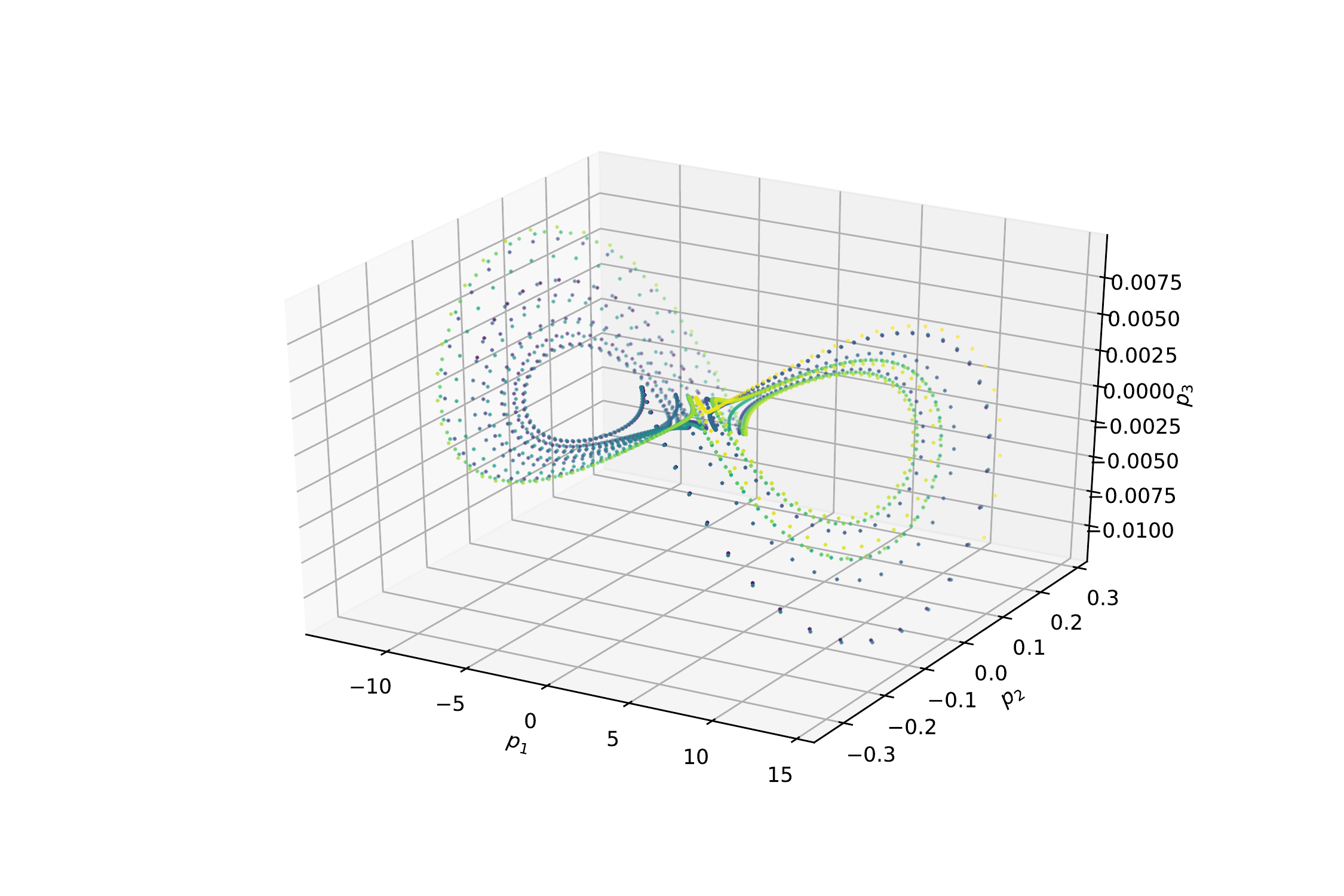}
 \caption{Projection onto the first three principal components of the reservoir states in the interval $(20,40)$ for a reservoir system driven by the first component of the Lorenz system. Points at the start of the trajectory are purple, and points later on are yellow.}
 \label{fig::PRS_lorenz_reservoir}
\end{figure}

\paragraph{The Lorenz system. Attractor reconstruction.} The Lorenz system with the parameter values given in the original paper \cite{lorenz1963deterministic} is determined by the differential equation
\begin{align*}
    \dot{u} &= 10(u-v), \\
    \dot{v} &= u(28-w) - v, \\
    \dot{w} &= uv - 8w/3.
\end{align*}
The discretization of this differential equation with time step $h$ yields a time evolution operator given by
\begin{align*}
    \phi(u_0,v_0,w_0) = (u_0,v_0,w_0) + \int_0^h (\dot{u}(t),\dot{v}(t),\dot{w}(t)) \ dt
\end{align*}
where the curve $(u(t),v(t),w(t))$ solves the Lorenz equations with initial condition $(u_0,v_0,w_0)$. In this paragraph we follow the same modeling prescriptions that we used for the R\"ossler and Van der Pol systems. We take the initial condition  as $ (u_0,v_0,w_0) = (0,1,1.05) $
and time step $h = 0.01$ and  the result for $T=40$ is plotted in Figure \ref{fig::PRS_lorenz}.
Regarding the reservoir embedding we use the same dimensionality that we took for R\"ossler and we also use the $u $ component as the input for state generation. A depiction of the projection of the corresponding states $x_t$ onto the first 3 principal components  is shown in Figure \ref{fig::PRS_lorenz_x}.
The corresponding states $\mathbf{x}_t $  in the interval $(20,40)$, that is, after a washout period of $(0,20)$ are plotted in Figure \ref{fig::PRS_lorenz_reservoir} after a projection onto the first three principal components.

\begin{figure}[h!]
\centering
\includegraphics[width=.75\textwidth]{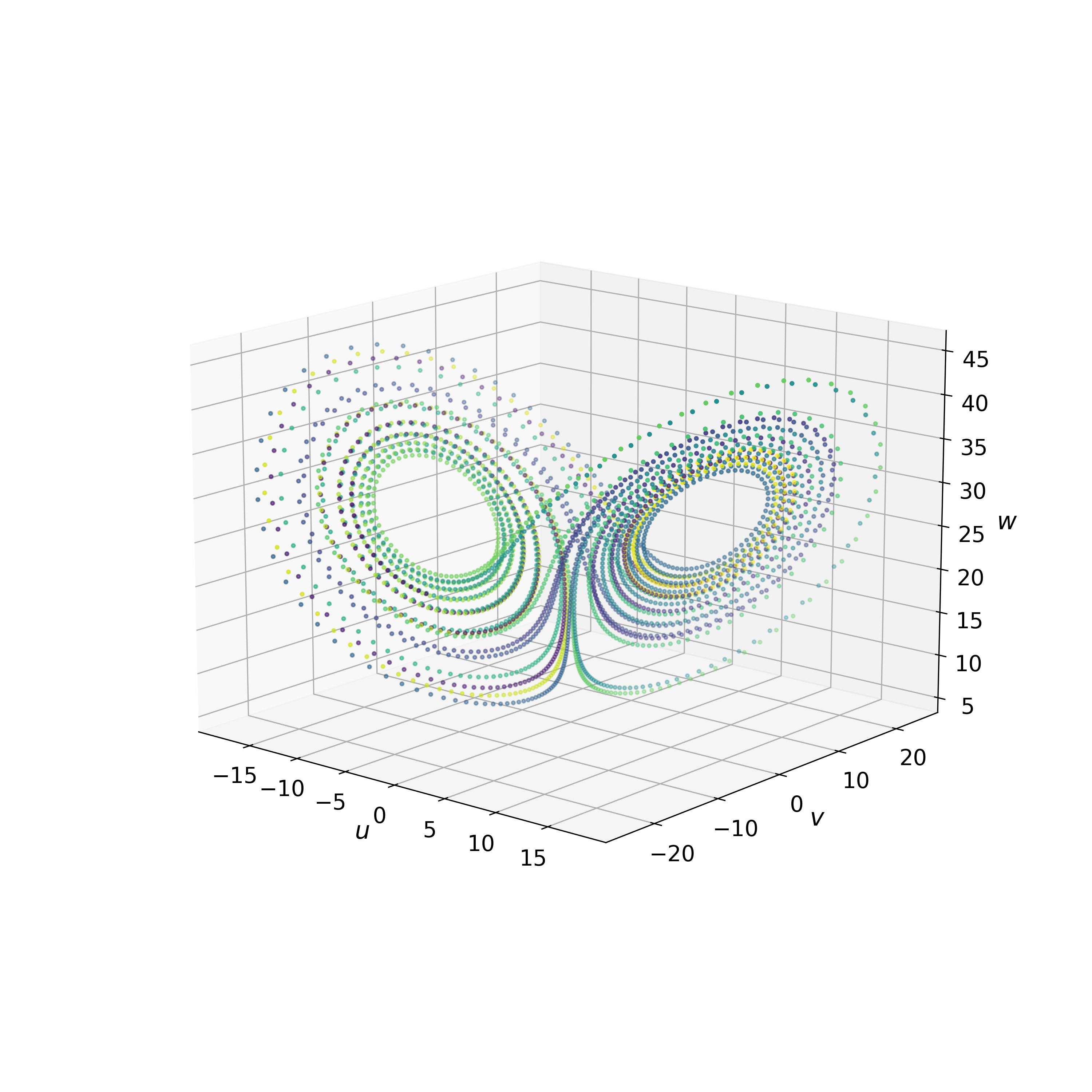}
\caption{A trajectory of the Lorenz system plotted for times in the interval $(1020,1050)$. Points at the start of the trajectory are purple, and points later on are yellow.}
\label{fig::Lorenz_traj_scatter}
\end{figure}

\paragraph{The Lorenz system. Forecasting in the presence of noise.} 
In this paragraph we follow closely the  experiment design that we previously used for attractor reconstruction but with a few modifications that we now list below. More precisely, we use the same parameters, initial point, and time step, but we consider $T=11000$ time units with a $1000$ time units long washout period. The result for this particularly choice of time interval is plotted in Figure~\ref{fig::Lorenz_traj_scatter}.
Consider now a reservoir system constructed according to the following prescription:
\begin{enumerate}
    \item Randomly generate a random orthogonal 20 by 20 matrix $A'$ drawn from the unique invariant Haar distribution in the Lie group ${\rm O}(20)$.
    \item Define the reservoir matrix $A := 0.9 \cdot A' / \lVert A' \rVert$.
    \item Randomly generate a vector $\mathbf{C'}\in \mathbb{R} ^{20}$  with IID uniform entries in the interval $[-1,1]$.
    \item Define the reservoir input vector $\mathbf{C} := \mathbf{C'}/ \lVert \mathbf{C'} \rVert$.
\end{enumerate}

\begin{figure}
  \centering
    \includegraphics[width=0.7\textwidth]{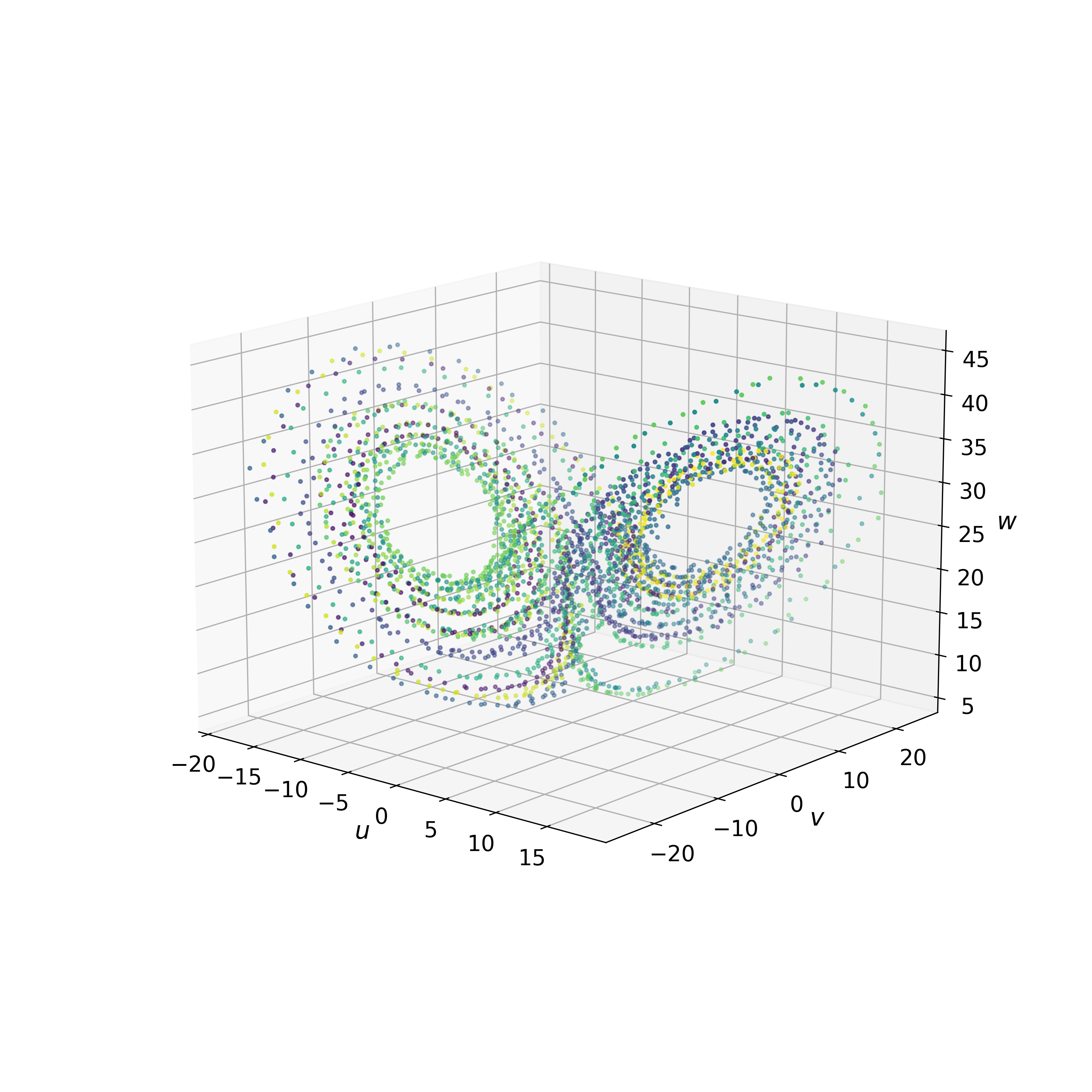}
 \caption{A trajectory of the Lorenz system with the $u$-component contaminated with  additive Gaussian noise and the original $v$ and $w$ components plotted for times in the interval $(1020,1050)$.}
 \label{fig::Lorenz_traj_ucompnoisy}
\end{figure}

We now choose a readout $h: \mathbb{R}^M \longrightarrow \mathbb{R}^m  $ for the state map using a deep neural network with $10$ hidden layers of $20$ neurons each and taking a scaled logistic map as activation function of the form  $\sigma(s)=(z_{{\rm max}}-z_{{\rm min}}) / (1 + e^{-s})+z_{{\rm min}}$. This readout is trained using states that are obtained with inputs of the form $u_t':=u_t+\epsilon_t$, where $u_t$ is the $u$-component of the Lorenz system and $\epsilon_t$ is a Gaussian distributed random variable with mean zero and variance $0.25$ for all $t \in (1000,11000)$. We illustrate a Lorenz trajectory for $(u',v,w)$ in Figure~\ref{fig::Lorenz_traj_ucompnoisy}.
 The corresponding states $\mathbf{x}_t $ are then subsequently used as the input of the deep neural network. The weights of the neural network are obtained by solving the empirical risk minimization problem where the one-step ahead shifted time series of the original $u$-component of the Lorenz system is taken as the target  and the mean squared error is taken as the empirical risk. The learning task hence consists in the  filtering of the noisy input and in the one-step ahead forecasting of the $u$-component of the dynamical Lorenz system. We implement the training of the deep neural network with the help of the Adam Optimizer in Keras in eight iterations with  early stopping parametrized by the patience parameter of 500 epochs. Each iteration consists of 7000 epochs of batch size 10000 and the learning rate is taken in a decreasing manner for each subsequent iteration out of the set of values $\{5{\rm e-}3,3{\rm e-}3,1{\rm e-}3,9{\rm e-}4,7{\rm e-}4,5{\rm e-}4,5{\rm e-}5,3{\rm e-}5\}$.

\begin{figure}
  \centering
    \includegraphics[width=0.55\textwidth]{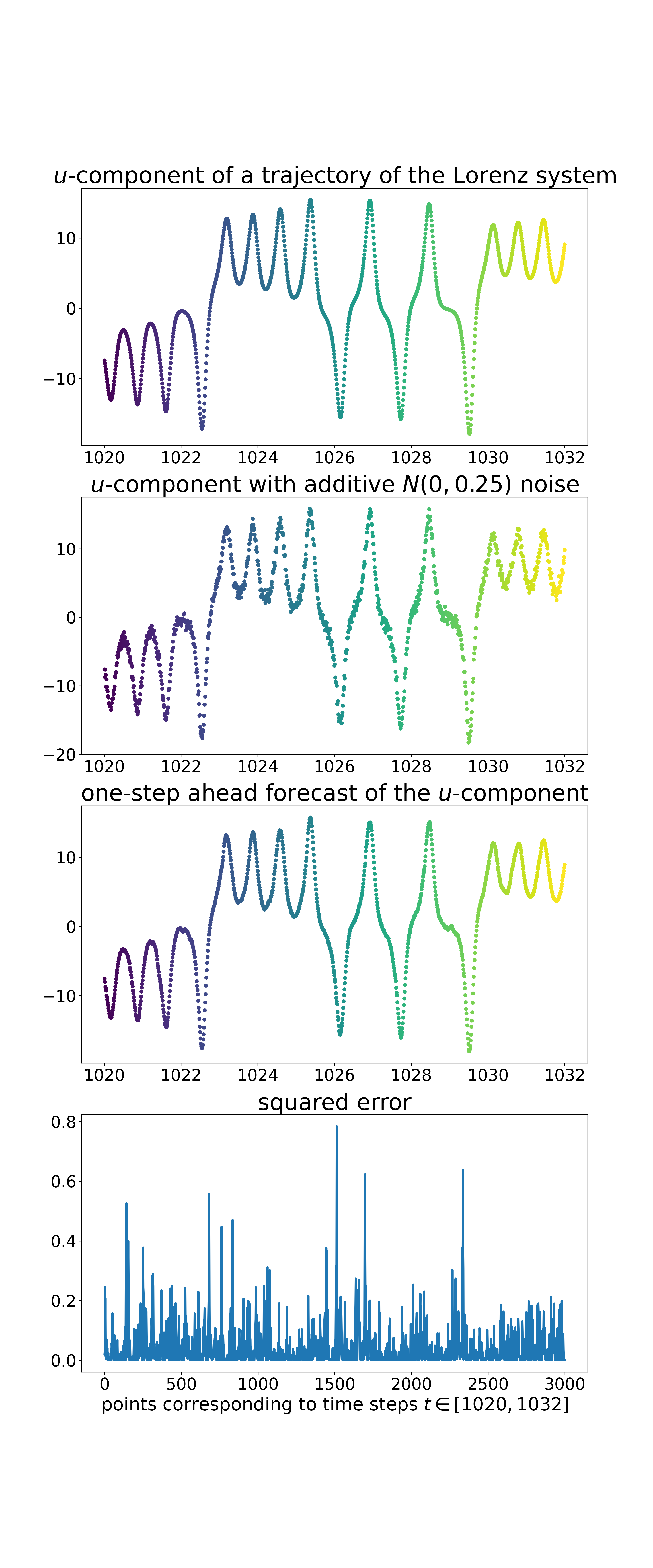}
 \caption{Filtered and forecasted $u$-component of the Lorenz system. The reservoir system is presented with a noisy version of the $u$-component time series as input and equipped with a trained deep neural network designed to filter and one-step-ahead forecast the time series.}
 \label{fig::Lorenz_traj_ucomponent}
\end{figure}

\begin{figure}[h!]
  \centering
    \includegraphics[width=0.7\textwidth]{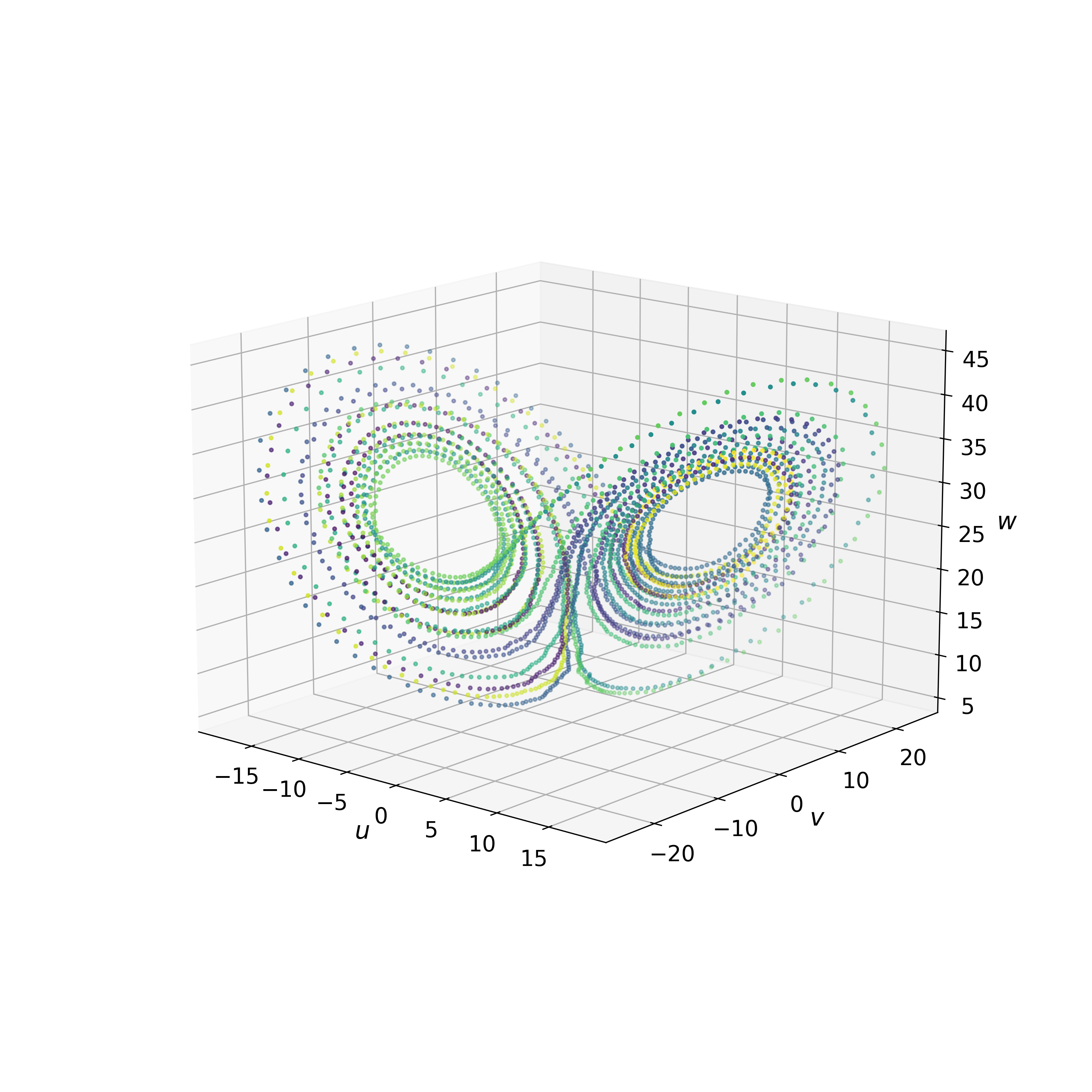}
 \caption{A trajectory of the Lorenz system with the filtered and forecasted $u$-component plotted for times in the interval $(1020,1050)$. The reservoir system is presented with a noisy version of the $u$-component time series as input and equipped with a trained deep neural network designed to filter and one-step-ahead forecast the time series.}
 \label{fig::Lorenz_traj_filtfrcst}
\end{figure}

The results on a testing sample are demostrated in Figure~\ref{fig::Lorenz_traj_ucomponent}. We complement the illustration with Figure~\ref{fig::Lorenz_traj_filtfrcst} which supports the pertinence of the methodology that we propose for denoising and forecasting since the reconstructed attractor is difficult to distinguish from the original one given in Figure~\ref{fig::Lorenz_traj_scatter}.

\section{Appendices}

\subsection{Elementary fact in linear algebra}

\begin{lemma}
\label{invertible perturbation}
Let $A$ and $B$ two square matrices of the same size such that $\det (A)=0 $ and $\det (B)\neq 0  $. Then, there exists $\varepsilon>0  $ such that 
\begin{equation*}
\det (A- \varepsilon B)\neq 0.
\end{equation*}
\end{lemma}

\noindent\textbf{Proof.\ \ } Consider the singular matrix $B ^{-1}A $ and let $\lambda_0  $ be its non-zero eigenvalue that has the smallest absolute value. Then, for any $0< \varepsilon< | \lambda_0| $ we necessarily have that $\det  \left(B ^{-1}A - \varepsilon \mathbb{I} \right)\neq 0 $ because otherwise $\varepsilon $ would be an eigenvalue of $B ^{-1}A $ which is impossible by the minimality of $\lambda_0 $. This implies that $C:= B ^{-1}A - \varepsilon \mathbb{I} $ is invertible and hence so is $BC=A- \varepsilon B $, as required. \quad $\blacksquare$

\medskip

As a corollary of this lemma we can conclude that if $\mathcal{V}:= \left\{\mathbf{v}_1, \ldots, \mathbf{v}_n\right\} $ and $\mathcal{W}:= \left\{\mathbf{w}_1, \ldots, \mathbf{w}_n\right\} $ are two $n$-sets of vectors in ${\Bbb R}^n $, then there exists $\varepsilon>0  $ such that the set $\left\{\mathbf{v}_1+ \varepsilon \mathbf{w} _1, \ldots, \mathbf{v}_n+ \varepsilon \mathbf{w} _n\right\} $ is made of linearly independent vectors. This fact is used at the end of the proof of Step 1 of Theorem \ref{Theorem immersion}.

\subsection{Two lemmas about random matrices}

\begin{lemma}
\label{first lemma polynomial}
Let $X_1, \ldots , X_n$ be independent real-valued non-singular random variables and let $p$ be a non-trivial polynomial in $n$ complex variables. Then
\begin{align*}
    \mathbb{P}\left(p(X_1, \ldots, X_n) = 0\right) = 0.
\end{align*}
\label{non_zero_polynomial_lemma}
\end{lemma}

\noindent\textbf{Proof.\ \ }Define $\mu_j(\cdot) := \mathbb{P}(X_j \in  \cdot),$ $j = 1, \ldots, n$, and let $Z = \{ \mathbf{x} \in \mathbb{C}^n \mid \ p(\mathbf{x}) = 0 \}$ be the set of complex roots of the polynomial $p$. Then, since  $X_1, \ldots , X_n$ are independent we have that
\begin{align*}
    \mathbb{P}\left(p(X_1, \ldots , X_n) = 0\right) = \mathbb{P}\left((X_1 , \ldots , X_n) \in Z\right) = (\mu_1 \otimes \ldots \otimes \mu_n)(Z).
\end{align*}
We now proceed by induction over $n$. For $n = 1$ we have that $\mathbb{P}\left(p(X_1) = 0\right) = \mu_1(Z) = 0$ since $Z$ is finite and $X_1$ is non singular. Let the claim be true for $n-1$. For fixed $x_1 \in \mathbb{C}$ set $p_{x_1}(x_2 , \ldots , x_n) := p(x_1 , \ldots , x_n)$ and let
\begin{align*}
    Z_{x_1} := \{ (x_2 , \ldots , x_n) \in \mathbb{C}^{n-1} \mid \ p_{x_1}(x_2 , \ldots , x_n) = 0\}.
\end{align*}
The set $F := \{ x_1 \in \mathbb{R} \mid \ p_{x_1} \equiv 0 \}$ is a finite set, so
\begin{equation*}
    \mathbb{P}\left(p(X_1, \ldots , X_n) = 0\right) = \int_{\mathbb{C}}(\mu_2 \otimes \ldots \otimes \mu_n)(Z_{x_1}) d \mu_1(x_1) = \int_{\mathbb{C} - F} (\mu_2 \otimes \ldots \otimes \mu_n)(Z_{x_1}) d \mu_1(x_1) = 0,
\end{equation*}
since we assumed that $X_1$ is non-singular and $(\mu_2 \otimes , \ldots , \otimes \mu_n)(Z_{x_1}) = 0$ for $x_1 \notin F$ by the induction hypothesis.
\quad $\blacksquare$

\begin{lemma}
\label{polynomial_lemma}
Let $N \in \mathbb{R}^N $, let  $A$ be a real $N \times N$ matrix, and let $\mathbf{C}$ be a random vector in $\mathbb{R}^N$.  Assume the entries of $A$ and $\mathbf{C}$ have been drawn using independent non-singular real valued random variables. Moreover, let $p_1 , \ldots , p_n \in \mathbb{C}[x]$ be linearly independent polynomials in one variable of degree at most $n-1$. Then
\begin{align*}
    \mathbb{P}(\det\left(p_1(A)\mathbf{C} | p_2(A)\mathbf{C} | \cdots | p_n(A)\mathbf{C}\right) = 0) = 0.
\end{align*}
Equivalently, the vectors $p_1(A)\mathbf{C} , p_2(A)\mathbf{C} , \ldots , p_n(A)\mathbf{C}$ are linearly independent almost surely.
\end{lemma}

\noindent\textbf{Proof.\ \ } The expression $\det\left(p_1(A)\mathbf{C} | \cdots | p_n(A)\mathbf{C}\right)$ is a polynomial $p$ in the $n^2 + n$ variables $a_{ij}$ and $b_j$, $i,j \in \left\{ 1 , \ldots , n\right\}$ that constitute the entries of $A$ and $\mathbf{C} $, respectively, and which in turn are by hypothesis non-singular random variables. As long as the polynomial $p$ is not identically zero, the result follows directly from the Lemma \ref{non_zero_polynomial_lemma}. So all that remains is to show that $p$ is not identically zero, that is, that there exist particular choices of $A$ and $\mathbf{C} $ such that $\det\left(p_1(A)\mathbf{C} | \cdots | p_n(A)\mathbf{C}\right)$ is non-zero. So, we choose $\mathbf{C} = (1 , \ldots , 1)^{\top}$ and $A = \text{diag}(a_1 , \ldots , a_n)$ with distinct real numbers $a_1 , \ldots , a_n$. We expand the polynomials $p_j$ in terms of their coefficients $\gamma_{j1} , \ldots , \gamma_{jn}$ so
\begin{align*}
    p_j(x) = \sum_{k = 0}^{n-1} \gamma_{jk}x^k.
\end{align*}
The vectors $\gamma_j = (\gamma_{j1}, \ldots , \gamma_{jn})^{\top}$, $j \in \left\{1, \ldots , n\right\}$, are by hypothesis linearly independent. We now show that with these choices, the vectors $p_1(A)\mathbf{C} , p_2(A)\mathbf{C} , \ldots , p_n(A)\mathbf{C}$ are linearly independent and hence $\det\left(p_1(A)\mathbf{C} | p_2(A)\mathbf{C} | \cdots | p_n(A)\mathbf{C}\right) \neq 0$, as required. Indeed, let $c_1 , \ldots , c_n \in \mathbb{R}$ and suppose that $\sum_{j=1}^n c_j p_j(A) \mathbf{C} ={\bf 0}$. Additionally, we can write
\begin{align*}
    \sum_{j=1}^n c_j p_j(A) \mathbf{C} = \sum_{j=1}^n c_j \sum_{k=0}^{n-1} \gamma_{kj} A^k \mathbf{C} = \sum_{k=0}^{n-1} \bigg( \sum^n_{j=1} c_j \gamma_{jk} \bigg) 
    \begin{bmatrix}
        a_1^k \\ a_2^k \\ \vdots \\ a_n^k
    \end{bmatrix}
    = V\mathbf{x}
\end{align*}
where 
\begin{equation*}
    V = 
\left(
\begin{array}{ccccc}
        1 & a_1 & a_1^2 & \ldots & a_1^{n-1} \\
        1 & a_2 & a_2^2 & \ldots & a_2^{n-1} \\
        \vdots & \vdots & \vdots & \ddots & \vdots \\
        1 & a_n & a_n^2 & \ldots & a_n^{n-1}
\end{array}
\right)
\quad \mbox{and} \quad
    \mathbf{x} = \sum^n_{j=1} c_j \gamma_{j}.
\end{equation*}
Since the diagonal entries of $A$ are all different and the determinant of the Vandermonde matrix $V$ is given by 
\begin{equation*}
\det (V)=\prod_{1\leq i<j\leq n} \left(a _i- a _j\right)
\end{equation*}
we can conclude that $V$ is invertible and hence the identity $\sum_{j=1}^n c_j p_j(A) \mathbf{C} =V \mathbf{x}={\bf 0}$ implies that $\mathbf{x}= {\bf 0} $. By the linear independence of the vectors $\gamma_1, \ldots, \gamma _n $ we have that $c _1, \ldots, c _n=0 $ necessarily.
It follows that the vectors $p_1(A)\mathbf{C} , p_2(A)\mathbf{C} , \ldots , p_n(A)\mathbf{C}$ are linearly independent, as required. 
\quad $\blacksquare$

\medskip

\noindent {\bf Acknowledgments.} We thank Friedrich Philipp for kindly communicating to us the proof of Proposition \ref{random_matrix_lemma}. We also thank  Henrik Brautmeier for his assistance with some of the numerical illustrations in the paper.   AH is supported by a scholarship from the EPSRC Centre for Doctoral Training in Statistical Applied Mathematics at Bath (SAMBa), project EP/L015684/1. JPO acknowledges partial financial support  coming from the Swiss National Science Foundation (grant number 200021\_175801/1).

\bibliographystyle{wmaainf}

\begin{thebibliography}{Gono~20b}

\bibitem[Abra~67]{abraham:robbin}
R.~Abraham and J.~Robbin.
\newblock {\it {Transversal Mappings and Flows}}.
\newblock W. A. Benjamin, Inc, 1967.

\bibitem[Abra~88]{mta}
R.~Abraham, J.~E. Marsden, and T.~S. Ratiu.
\newblock {\it {Manifolds, Tensor Analysis, and Applications}}.
\newblock Vol.~75, Applied Mathematical Sciences. Springer-Verlag, 1988.

\bibitem[Apos~74]{Apostol:analysis}
T.~Apostol.
\newblock {\it {Mathematical Analysis}}.
\newblock Addison Wesley, second Ed., 1974.

\bibitem[Bocc~02]{boccaletti:reports:2002}
S.~Boccaletti, J.~Kurths, G.~Osipov, D.~L. Valladares, and C.~S. Zhou.
\newblock ``{The synchronization of chaotic systems}''.
\newblock {\it Physics Reports}, Vol.~366, pp.~1--101, 2002.

\bibitem[Boot~03]{boothby2003introduction}
W.~M. Boothby.
\newblock {\it {An Introduction to Differentiable Manifolds and Riemannian
  Geometry}}.
\newblock Academic Press, Inc., second rev Ed., 2003.

\bibitem[Carm~92]{do:carmo:1993}
M.~P. do~Carmo.
\newblock {\it {Riemannian Geometry}}.
\newblock Birkh{\"{a}}user Boston, 1992.

\bibitem[Carr~18]{carroll2018using}
T.~L. Carroll.
\newblock ``{Using reservoir computers to distinguish chaotic signals}''.
\newblock {\it Physical Review E}, Vol.~98, No.~5, p.~52209, 2018.

\bibitem[Erog~17]{eroglu2017synchronisation}
D.~Eroglu, J.~S.~W. Lamb, and T.~Pereira.
\newblock ``{Synchronisation of chaos and its applications}''.
\newblock {\it Contemporary Physics}, Vol.~58, No.~3, pp.~207--243, 2017.

\bibitem[Gaut~21]{gauthier2021next}
D.~J. Gauthier, E.~Bollt, A.~Griffith, and W.~A.~S. Barbosa.
\newblock ``{Next Generation Reservoir Computing}''.
\newblock {\it arXiv preprint arXiv:2106.07688}, 2021.

\bibitem[Gono~20a]{RC15}
L.~Gonon, L.~Grigoryeva, and J.-P. Ortega.
\newblock ``{Memory and forecasting capacities of nonlinear recurrent
  networks}''.
\newblock {\it Physica D}, Vol.~414, No.~132721, pp.~1--13., 2020.

\bibitem[Gono~20b]{RC8}
L.~Gonon and J.-P. Ortega.
\newblock ``{Reservoir computing universality with stochastic inputs}''.
\newblock {\it IEEE Transactions on Neural Networks and Learning Systems},
  Vol.~31, No.~1, pp.~100--112, 2020.

\bibitem[Gono~21]{RC20}
L.~Gonon and J.-P. Ortega.
\newblock ``{Fading memory echo state networks are universal}''.
\newblock {\it Neural Networks}, Vol.~138, pp.~10--13, 2021.

\bibitem[Grig~18]{RC7}
L.~Grigoryeva and J.-P. Ortega.
\newblock ``{Echo state networks are universal}''.
\newblock {\it Neural Networks}, Vol.~108, pp.~495--508, 2018.

\bibitem[Grig~20a]{RC18}
L.~Grigoryeva, A.~G. Hart, and J.-P. Ortega.
\newblock ``{Chaos on compact manifolds: Differentiable synchronizations beyond
  Takens}''.
\newblock {\it Preprint arXiv:2010.03218}, 2020.

\bibitem[Grig~20b]{RC16}
L.~Grigoryeva and J.-P. Ortega.
\newblock ``{Dimension reduction in recurrent networks by canonicalization}''.
\newblock {\it Preprint arXiv:2007.12141}, 2020.

\bibitem[Hart~20]{hart:ESNs}
A.~G. Hart, J.~L. Hook, and J.~H.~P. Dawes.
\newblock ``{Embedding and approximation theorems for echo state networks}''.
\newblock {\it Neural Networks}, Vol.~128, pp.~234--247, 2020.

\bibitem[Hart~21]{allen:tikhonov}
A.~G. Hart, J.~L. Hook, and J.~H.~P. Dawes.
\newblock ``{Echo State Networks trained by Tikhonov least squares are
  L2($\mu$) approximators of ergodic dynamical systems}''.
\newblock {\it Physica D: Nonlinear Phenomena}, p.~132882, 2021.

\bibitem[Hirs~76]{Hirsch:book}
M.~W. Hirsch.
\newblock {\it {Differential Topology}}.
\newblock Springer Verlag, 1976.

\bibitem[Huke~06]{huke:2006}
J.~P. Huke.
\newblock ``{Embedding nonlinear dynamical systems: a guide to Takens'
  theorem}''.
\newblock Tech. Rep., Manchester Institute for Mathematical Sciences. The
  University of Manchester, 2006.

\bibitem[Jaeg~04]{Jaeger04}
H.~Jaeger and H.~Haas.
\newblock ``{Harnessing Nonlinearity: Predicting Chaotic Systems and Saving
  Energy in Wireless Communication}''.
\newblock {\it Science}, Vol.~304, No.~5667, pp.~78--80, 2004.

\bibitem[Jaeg~10]{jaeger2001}
H.~Jaeger.
\newblock ``{The `echo state' approach to analysing and training recurrent
  neural networks with an erratum note}''.
\newblock Tech. Rep., German National Research Center for Information
  Technology, 2010.

\bibitem[Kalm~10]{Kalman2010}
R.~Kalman.
\newblock ``{Lectures on Controllability and Observability}''.
\newblock In: {\it Controllability and Observability}, pp.~1--149, Springer
  Berlin Heidelberg, Berlin, Heidelberg, 2010.

\bibitem[Kant~03]{kantz:Sreiber}
H.~Kantz and T.~Schreiber.
\newblock {\it {Nonlinear Time Series Analysis}}.
\newblock Cambridge University Press, second Ed., 2003.

\bibitem[Kryl~31]{krylov1931numerical}
A.~N. Krylov.
\newblock ``{On the numerical solution of equation by which are determined in
  technical problems the frequencies of small vibrations of material
  systems}''.
\newblock {\it News Acad. Sci. USSR}, Vol.~7, pp.~491--539, 1931.

\bibitem[Kupk~63]{kupka1963contributiona}
I.~Kupka.
\newblock ``{Contributiona la th{\'{e}}orie des champs
  g{\'{e}}n{\'{e}}riques}''.
\newblock {\it Contributions to differential equations}, Vol.~2, pp.~457--484,
  1963.

\bibitem[Lax~02]{lax:functional:analysis}
P.~Lax.
\newblock {\it {Functional Analysis}}.
\newblock Wiley-Interscience, 2002.

\bibitem[Lore~63]{lorenz1963deterministic}
E.~N. Lorenz.
\newblock ``{Deterministic nonperiodic flow}''.
\newblock 1963.

\bibitem[Lu~18]{Ott2018}
Z.~Lu, B.~R. Hunt, and E.~Ott.
\newblock ``{Attractor reconstruction by machine learning}''.
\newblock {\it Chaos}, Vol.~28, No.~6, 2018.

\bibitem[Lu~20]{lu:bassett:2020}
Z.~Lu and D.~S. Bassett.
\newblock ``{Invertible generalized synchronization: A putative mechanism for
  implicit learning in neural systems}''.
\newblock {\it Chaos}, Vol.~30, No.~063133, 2020.

\bibitem[Luko~09]{lukosevicius}
M.~Luko{\v{s}}evi{\v{c}}ius and H.~Jaeger.
\newblock ``{Reservoir computing approaches to recurrent neural network
  training}''.
\newblock {\it Computer Science Review}, Vol.~3, No.~3, pp.~127--149, 2009.

\bibitem[Maas~00]{Maass2000}
W.~Maass and E.~D. Sontag.
\newblock ``{Neural Systems as Nonlinear Filters}''.
\newblock {\it Neural Computation}, Vol.~12, No.~8, pp.~1743--1772, aug 2000.

\bibitem[Maas~02]{maass1}
W.~Maass, T.~Natschl{\"{a}}ger, and H.~Markram.
\newblock ``{Real-time computing without stable states: a new framework for
  neural computation based on perturbations}''.
\newblock {\it Neural Computation}, Vol.~14, pp.~2531--2560, 2002.

\bibitem[Maas~04]{corticalMaass}
W.~Maass, T.~Natschl{\"{a}}ger, and H.~Markram.
\newblock ``{Fading memory and kernel properties of generic cortical
  microcircuit models}''.
\newblock {\it Journal of Physiology Paris}, Vol.~98, No.~4-6 SPEC. ISS.,
  pp.~315--330, 2004.

\bibitem[Maas~07]{MaassUniversality}
W.~Maass, P.~Joshi, and E.~D. Sontag.
\newblock ``{Computational aspects of feedback in neural circuits}''.
\newblock {\it PLoS Computational Biology}, Vol.~3, No.~1, p.~e165, 2007.

\bibitem[Manj~13]{Manjunath:Jaeger}
G.~Manjunath and H.~Jaeger.
\newblock ``{Echo state property linked to an input: exploring a fundamental
  characteristic of recurrent neural networks}''.
\newblock {\it Neural Computation}, Vol.~25, No.~3, pp.~671--696, 2013.

\bibitem[Manj~20]{manjunath:prsl}
G.~Manjunath.
\newblock ``{Stability and memory-loss go hand-in-hand: three results in
  dynamics {\&} computation}''.
\newblock {\it To appear in Proceedings of the Royal Society London Ser. A
  Math. Phys. Eng. Sci.}, pp.~1--25, 2020.

\bibitem[Matt~92]{Matthews:thesis}
M.~B. Matthews.
\newblock {\it {On the Uniform Approximation of Nonlinear Discrete-Time
  Fading-Memory Systems Using Neural Network Models}}.
\newblock PhD thesis, ETH Z{\"{u}}rich, 1992.

\bibitem[Matt~93]{Matthews1993}
M.~B. Matthews.
\newblock ``{Approximating nonlinear fading-memory operators using neural
  network models}''.
\newblock {\it Circuits, Systems, and Signal Processing}, Vol.~12, No.~2,
  pp.~279--307, jun 1993.

\bibitem[Munk~14]{Munkres:topology}
J.~Munkres.
\newblock {\it {Topology}}.
\newblock Pearson, second Ed., 2014.

\bibitem[Nats~02]{Natschlager:117806}
T.~Natschl{\"{a}}ger, W.~Maass, and H.~Markram.
\newblock ``{The "Liquid Computer": a novel strategy for real-time computing on
  time series}''.
\newblock {\it Special Issue on Foundations of Information Processing of
  TELEMATIK}, Vol.~8, No.~1, pp.~39--43, 2002.

\bibitem[Ott~02]{ott2002chaos}
E.~Ott.
\newblock {\it {Chaos in Dynamical Systems}}.
\newblock Cambridge University Press, second Ed., 2002.

\bibitem[Path~17]{pathak:chaos}
J.~Pathak, Z.~Lu, B.~R. Hunt, M.~Girvan, and E.~Ott.
\newblock ``{Using machine learning to replicate chaotic attractors and
  calculate Lyapunov exponents from data}''.
\newblock {\it Chaos}, Vol.~27, No.~12, 2017.

\bibitem[Path~18]{Pathak:PRL}
J.~Pathak, B.~Hunt, M.~Girvan, Z.~Lu, and E.~Ott.
\newblock ``{Model-Free Prediction of Large Spatiotemporally Chaotic Systems
  from Data: A Reservoir Computing Approach}''.
\newblock {\it Physical Review Letters}, Vol.~120, No.~2, p.~24102, 2018.

\bibitem[Peco~97]{pecora:synch}
L.~M. Pecora, T.~L. Carroll, G.~A. Johnson, D.~J. Mar, and J.~F. Heagy.
\newblock ``{Fundamentals of synchronization in chaotic systems, concepts, and
  applications}''.
\newblock {\it Chaos}, Vol.~7, No.~4, pp.~520--543, 1997.

\bibitem[Rulk~95]{rulkov1995generalized}
N.~F. Rulkov, M.~M. Sushchik, L.~S. Tsimring, and H.~D.~I. Abarbanel.
\newblock ``{Generalized synchronization of chaos in directionally coupled
  chaotic systems}''.
\newblock {\it Physical Review E}, Vol.~51, No.~2, p.~980, 1995.

\bibitem[Saue~91]{sauer1991embedology}
T.~Sauer, J.~A. Yorke, and M.~Casdagli.
\newblock ``{Embedology}''.
\newblock {\it Journal of Statistical Physics}, Vol.~65, No.~3, pp.~579--616,
  1991.

\bibitem[Smal~63]{smale1963stable}
S.~Smale.
\newblock ``{Stable manifolds for differential equations and
  diffeomorphisms}''.
\newblock {\it Annali della Scuola Normale Superiore di Pisa-Classe di
  Scienze}, Vol.~17, No.~1-2, pp.~97--116, 1963.

\bibitem[Sont~98]{sontag:book}
E.~Sontag.
\newblock {\it {Mathematical Control Theory: Deterministic Finite Dimensional
  Systems}}.
\newblock Springer-Verlag, 1998.

\bibitem[Take~81]{takensembedding}
F.~Takens.
\newblock ``{Detecting strange attractors in turbulence}''.
\newblock pp.~366--381, Springer Berlin Heidelberg, 1981.

\bibitem[Tana~19]{tanaka:review}
G.~Tanaka, T.~Yamane, J.~B. H{\'{e}}roux, R.~Nakane, N.~Kanazawa, S.~Takeda,
  H.~Numata, D.~Nakano, and A.~Hirose.
\newblock ``{Recent advances in physical reservoir computing: A review}''.
\newblock {\it Neural Networks}, Vol.~115, pp.~100--123, 2019.

\bibitem[Verz~20]{Verzelli2020b}
P.~Verzelli, C.~Alippi, and L.~Livi.
\newblock ``{Learn to Synchronize, Synchronize to Learn}''.
\newblock 2020.

\end{thebibliography}

\end{document}